\documentclass{amsart}
\usepackage{amsthm}
\usepackage{amsmath,amssymb}
\usepackage{url,float}
\usepackage{xcolor}
\usepackage{graphicx}
\usepackage{enumerate}
\usepackage{algorithm}
\usepackage{algpseudocode}
\usepackage{xfrac}
\usepackage[colorlinks]{hyperref}
\usepackage{tabularx}
\usepackage{tikz-cd}
\usepackage{tikz}
\usepackage{enumitem}
\usepackage{pifont}
\usepackage{diagbox}
\usepackage{multirow}
\usepackage{makecell}
\usepackage[most]{tcolorbox}
\usepackage{longtable}
\usepackage{pdflscape}
\usepackage{geometry}
\usepackage{subcaption}
\usepackage{url}

\newcommand{\F}{\mathbb{F}}
\newtheorem{thm}{Theorem}[section]
\newtheorem{cor}[thm]{Corollary}
\newtheorem{lem}[thm]{Lemma}
\newtheorem{rem}[thm]{Remark}
\newtheorem{defn}[thm]{Definition}
\newtheorem{pro}[thm]{Proposition}

\newtheorem{example}[thm]{Example}
\newtheorem{conj}[thm]{Conjecture}
\newtheorem{question}[thm]{Question}

\newcommand{\Z}{\mathbb{Z}}

\newcounter{alg}

\title[On the Lattice Tilings by Lee Balls]{On the Lattice Tilings by Lee Balls: Q-Polynomials, Group-Ring Equations, and a Decision Algorithm}
\author[A.\ Xiao]{Ang Xiao}
\address{College of Science, National University of Defense Technology, 410073 Changsha, China}
\email{xiaoang\_nudt@163.com}
\author[Y.\ Zhou]{Yue Zhou}
\address{College of Science, National University of Defense Technology, 410073 Changsha, China}
\email{yue.zhou.ovgu@gmail.com}
\urladdr{https://sites.google.com/site/yuejoezhou/}

\begin{document}
	\begin{abstract}
		More than 50 years ago, Golomb and Welch conjectured that no tilings of $\mathbb{Z}^n$ by Lee balls of radius $r$ exist for $n \geq 3$ and $r \geq 2$. This conjecture remains open even under the additional assumption that the tiling is a lattice tiling. In general, there are two known effective approaches to the lattice case which are based on symmetric polynomials and on the arithmetic of group-ring equations, respectively.

		In this paper, we develop both approaches. First, we determine all coefficients of the Q-polynomials associated with Lee balls; previously, only leading coefficients and certain special cases were known. Second, for arbitrary radii, we give necessary and sufficient conditions for the existence of lattice tilings by Lee balls in terms of group-ring equations, extending earlier results for radii $r=2$ and $3$.

		We then derive new necessary conditions for such tilings. For every radius $3\le r\le 50$, we prove that the set of dimensions for which no lattice tiling exists has positive density in \(\mathbb{Z}^+\). We also give an algorithm that combines the available necessary conditions and apply it to all pairs $(n,r)$ for which the corresponding Lee ball of radius $r$ and dimension $n$ has size less than $2^{30}$. Among these $1{,}358$ pairs, the combined method excludes $1{,}356$, leaving only $(n,r)=(35,3)$ and $(805,3)$.
	\end{abstract}
	\maketitle
	\tableofcontents
	\section{Introduction}
	The study of packings and tilings in continuous and discrete spaces is a central topic in discrete geometry, with deep connections to number theory, group theory, and coding theory. A prominent packing problem asks for the maximum density of non-overlapping spheres in Euclidean space $\mathbb{R}^n$. In dimension $n=3$, Kepler's conjecture on maximum density was famously proven by Hales \cite{hales2005proof}. More recently, spectacular breakthroughs resolved higher-dimensional cases: Viazovska \cite{viazovska2017sphere} solved the problem in $n=8$ using modular forms, and Cohn, Kumar, Miller, Radchenko, and Viazovska \cite{cohn2017sphere} extended this to $n=24$ via the Leech lattice.
	
	When packing density reaches unity, the configuration becomes a tiling. Tiling problems often exhibit rigid algebraic structures. For example, Minkowski (1907) conjectured that in any lattice tiling of $\mathbb{R}^n$ by unit cubes, two cubes must share a complete $(n-1)$-dimensional facet \cite{minkowski1907diophantische}. Haj\'os \cite{hajos1941uber} affirmatively resolved this in 1941 by elegantly translating the geometric setup into a purely algebraic problem concerning the factorization of finite abelian groups.
	
	However, structural rigidity frequently breaks down as the dimension $n$ grows. This is illustrated by Keller's conjecture (1930) \cite{keller1930uber}, which dropped the lattice condition from Minkowski's proposition, conjecturing that \textit{any} tiling of $\mathbb{R}^n$ by identical cubes has two sharing a complete $(n-1)$-dimensional facet. Surprisingly, this fails in high dimensions. The breakthrough came from Corr\'adi and Szab\'o \cite{corradi1990combinatorial}, who reduced the continuous problem to finding a clique of size $2^n$ in the discrete Keller graph. Using this framework, Lagarias and Shor \cite{lagarias1992keller} found counterexamples for $n \ge 10$, which Mackey \cite{mackey2002cube} lowered to $n \ge 8$. Brakensiek et al.\ \cite{brakensiek2020resolution} completely resolved the problem, showing the conjecture holds only for $n \le 7$.
	
	Such high-dimensional phenomena persist in the integer lattice $\mathbb{Z}^n$. The study of translational tilings of $\mathbb{Z}^n$ by finite subsets (polyominoes) has revealed striking structural complexities. A central question was the periodic tiling conjecture: if a finite subset tiles $\mathbb{Z}^n$ by translations, must it admit a periodic tiling? This was recently shattered in full generality by Greenfeld and Tao \cite{greenfeld2024counterexample}. They explicitly constructed a single finite subset of $\mathbb{Z}^n$ (for sufficiently large $n$) that successfully tiles the lattice but remains strictly aperiodic. This landmark result illustrates that high-dimensional translational tilings can avoid regular, group-theoretic structures entirely.
	\subsection{Golomb--Welch Conjecture} Against this backdrop of high-dimensional irregularity and aperiodicity in $\mathbb{Z}^n$, it is natural to ask whether imposing strict symmetries on a discrete shape can enforce regular tiling behaviors. The most natural highly symmetric discrete analogue of the Euclidean ball (and the continuous cross-polytope) is the Lee ball.

	Denote by $\Z$ the ring of integers. For any two points $x=(x_1,\dots,x_n)$ and $y=(y_1,\dots,y_n)$ in $\Z^n$, the \emph{Lee metric} (also known as Manhattan or $\ell_1$ metric) in $\Z^n$ is given by
	$$d_L(x,y):=\sum_{i=1}^{n}|x_i-y_i|.$$
	The \emph{Lee ball} in $\Z^n$ centered at $\textbf{0}$ with radius $r$ is defined as
	$$B^n(r):=\{\textbf{x}\in\Z^n:d_L(\textbf{x},\textbf{0})\le r\}.$$
	If there exists a subset $C\subseteq\Z^n$ such that the collection $\mathcal{T}:=\bigsqcup_{c\in C}\{b+c:b\in B^n(r)\}$ forms a partition of $\Z^n$, then $\mathcal{T}$ is called a \emph{tiling} of $\Z^n$ by the Lee ball $B^n(r)$. Furthermore, if $C$ is a lattice, we call $\mathcal T$ a \emph{lattice tiling} (by Lee balls).

	The literature often refers to these configurations as ``sphere packings'' or ``tilings by Lee spheres,'' although the sets defined by the inequality above are balls rather than boundaries. Because our arguments distinguish a ball from its boundary, we use \emph{Lee ball} for the solid set and reserve \emph{Lee sphere} for the boundary $\{\mathbf{x}\in\Z^n:d_L(\mathbf{x},\mathbf{0})=r\}$.

	In \cite{golomb_perfect_1970}, Golomb and Welch initiated the systematic study of perfect Lee codes. They observed that translational tilings of $\mathbb{Z}^n$ by Lee balls exist under the assumption that $n=1,2$ or $r=1$. In the same paper, they also established the nonexistence of such tilings for sufficiently large radii (with a threshold $r_n$ depending on $n$, not explicitly specified), and conjectured that a tiling of $\Z^n$ by $B^n(r)$ exists if and only if $r=1$ or $n=1,2$.

	\begin{conj}[\cite{golomb_perfect_1970}, Golomb and Welch]
		There is no tiling of $ \Z^n $ by the Lee ball $ B^n(r) $ for $ n\ge3 $ and $ r\ge2$.
	\end{conj}

	The lattice restriction of the Golomb--Welch conjecture also has a direct interpretation in the degree--diameter problem for Abelian Cayley graphs. Let $G$ be a finite abelian group and let
	\[
	S=\{\pm g_1,\dots,\pm g_n\}
	\]
	be an inverse-closed generating set of size $2n$. If the Cayley graph $\operatorname{Cay}(G,S)$ has diameter at most $r$, then the homomorphism
	\[
	\varphi:\mathbb Z^n\longrightarrow G,
	\qquad
	(a_1,\dots,a_n)\longmapsto \sum_{i=1}^n a_i g_i,
	\]
	is surjective on $B^n(r)$. Consequently,
	\[
	|G|\le |B^n(r)|
	=\sum_{i=0}^{\min\{n,r\}}2^i\binom ni\binom ri.
	\]
	The right-hand side is the \emph{Abelian Cayley--Moore bound} for degree $2n$ and diameter $r$ \cite{dougherty_degree-diameter_2004}. Equality holds precisely when $\varphi$ is bijective on $B^n(r)$; by the standard homomorphism criterion recalled in Corollary~\ref{cor:P2-group_latticetiling}, this is equivalent to the existence of a lattice tiling of $\mathbb Z^n$ by $B^n(r)$ \cite{horak_diameter_2012}. Thus, for $n\ge3$ and $r\ge2$, existence of a lattice tiling is equivalent to attainment of the Abelian Cayley--Moore bound. For more degree--diameter problems in graph theory, we refer to the survey \cite{miller_siran_moore_2013}. We nevertheless keep lattice tilings as the main theme and use their language throughout the proofs.

	Let $\Z_q$ denote the ring of integers modulo $q$. A similar conjecture for tilings of $\Z_q^n$ by Lee balls was proposed in \cite[Section 3]{golomb_perfect_1970}. In fact, for any $q\geq 2r+1$, a tiling of $\Z_q^n$ by Lee balls of radius $r$ lifts to a $q$-periodic tiling of $\Z^n$. In the literature, the conjecture for $\Z_q^n$ is usually called the weak version and the one for $\Z^n$ is called the strong version.

\begin{figure}[htp]
		\centering
		\includegraphics[width=0.9\linewidth]{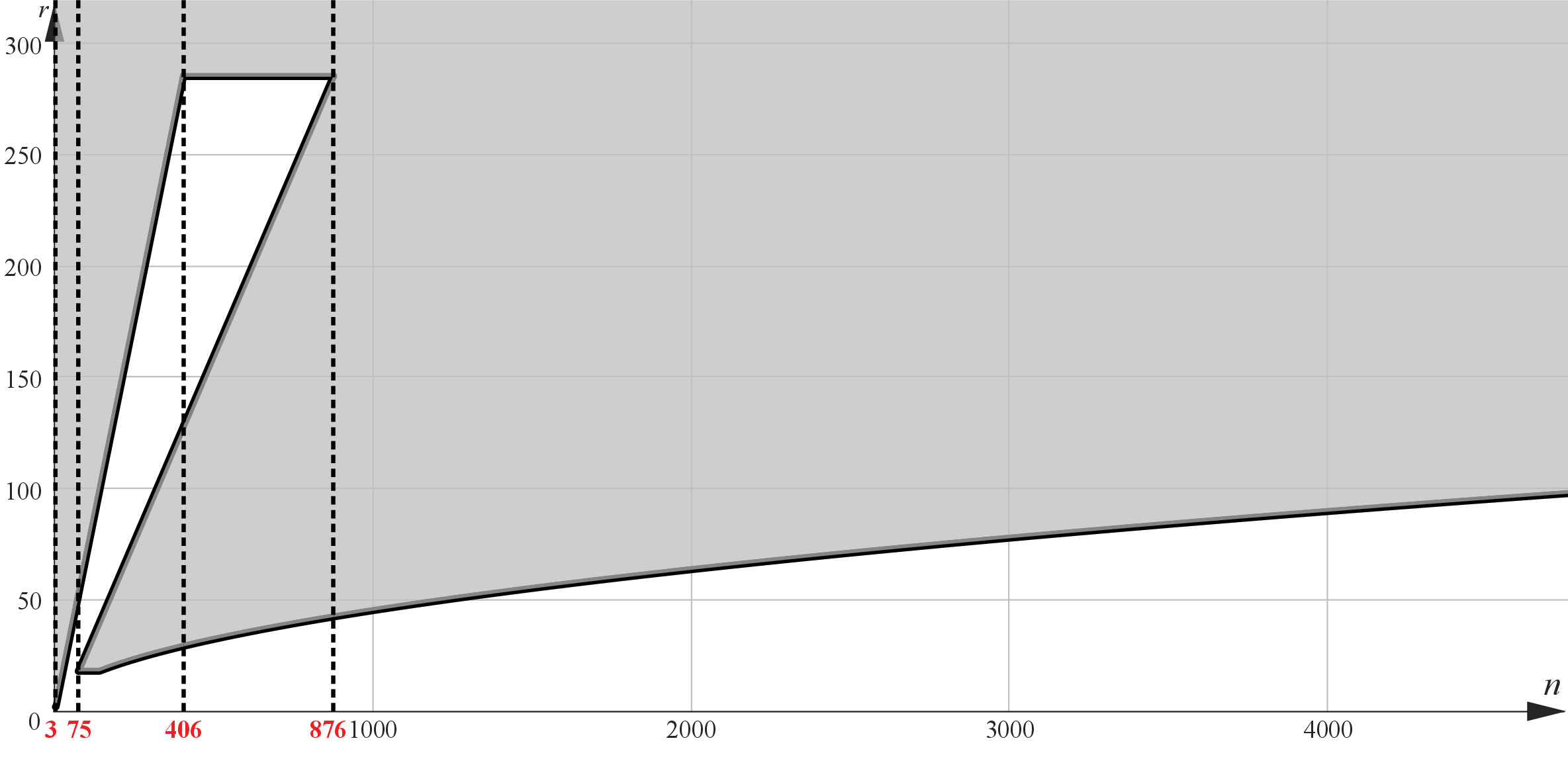}
		\caption{Proven and Open Cases of the Golomb--Welch Conjecture Summarized in Theorem \ref{thm:P1-HK}}
		\label{fig:P1-horak_kim}
	\end{figure}

	Both the weak and strong Golomb--Welch conjectures remain open. Much of the literature concerns bounds on $r_n$.
	In \cite{post_nonexistence_1975}, Post showed that there is no tiling of $\Z_q^n$ by $B^n(r)$ when $r\ge \frac{\sqrt{2}}{2}n-\frac{3}{4}\sqrt{2}-\frac{1}{2}$ for $n\ge6$, or when $r\ge n-1$ for $3\le n\le5$. This result was improved by Astola \cite{astola_elias-type_1982}, and later by Lepist\"o \cite{lepisto_modification_1981} who proved that tiling of $\Z_q^n$ by $B^n(r)$ with $q\geq 2r+1$ must satisfy $n\ge\frac{(r+2)^2}{2.1}$, where $r\ge285$. In \cite{kim_50_2018}, Horak and Kim dropped the periodicity condition on these bounds, and proved that there is no tiling of $\Z^n$ by Lee ball $B^n(r)$ for $r\ge18$ and $3r+21\le n\le\frac{1}{2}r^2-20$. The following theorem combines these bounds on $r_n$.
	\begin{thm}[\cite{kim_50_2018}, Horak and Kim]\label{thm:P1-HK}
		There is no tiling of $\Z^n$ by Lee ball $B^n(r)$ for
		\begin{itemize}
			\item $3 \leq n \leq 74$ and $\max\left\{ \frac{\sqrt{2}}{2}n - \frac{3}{4}\sqrt{2} - \frac{1}{2},2 \right\} \leq r$.
			\item $75 \leq n \leq 405$ and $\max\left\{ 18,\sqrt{2n + 40} \right\} \leq r \leq \frac{n - 21}{3}$ or $\frac{\sqrt{2}}{2}n - \frac{3}{4}\sqrt{2} - \frac{1}{2} \leq r$.
			\item $406 \leq n \leq 876$ and $\sqrt{2n + 40} \leq r \leq \frac{n - 21}{3}$ or $285 \leq r.$
			\item $n \geq 876$ and $\sqrt{2n + 40} \leq r$.
		\end{itemize}
	\end{thm}
	Based on Theorem \ref{thm:P1-HK}, we plot the parameter pairs $(n,r)$ for which there is no tiling of $\Z^n$ by $B^n(r)$, shown in gray in Fig.~\ref{fig:P1-horak_kim}, with vertical dashed lines at $n=3,75,406,876$.

	The nonexistence results summarized in Theorem~\ref{thm:P1-HK} are not exhaustive. For several additional values of $n$ and $r$, the Golomb--Welch conjecture has also been proved:
	\begin{itemize}
		\item $n=6$ and $r=2$; see \cite{horak_perfect_2009};
		\item $3\leq n \leq 5$ and $r\geq 2$; see \cite{gravier_non-existence_1998,horak_tilings_2009,horak_new_2014,spacapan_nonexistence_2007}.
	\end{itemize}
	Moreover, when $|B^n(r)|$ is prime, a result of Szegedy \cite{szegedy_algorithms_1998} implies that nonexistence of a lattice tiling also yields nonexistence of an arbitrary translational tiling for the same $n$ and $r$.

	\subsection{Lattice Tilings of $\Z^n$ by Lee Balls} In contrast to the results listed in Theorem \ref{thm:P1-HK}, the lattice cases of the Golomb--Welch conjecture are usually proved for a fixed radius $r$ while allowing the dimension $n$ to vary.

	There are two principal approaches to the lattice-tiling case: one is based on symmetric polynomials over finite fields, and the other is based on the arithmetic of group-ring equations.

	In \cite{kim_nonexistence_2017}, by using symmetric functions, Kim established a necessary condition for a lattice tiling of $\Z^n$ by Lee balls of radius $2$ that can be verified computationally. For instance, it leads to the proof of the Golomb--Welch conjecture for $r=2$ and $222$ particular values of $n$ for $3\leq n\leq 10^3$.
	For $r=3,4$, new criteria were obtained by Zhang and Ge in \cite{zhang_perfect_2017}.
	In \cite{qureshi_non-existence_2020}, Qureshi introduced a more general framework to extend the previous approaches. For any subset $ B\subseteq\Z^n $, define the Q-polynomial
	\begin{align*}
		Q_B^k(X):=\sum_{b\in B}\left\langle X,b\right\rangle^{2k},
	\end{align*}
	where $X=(X_1,\dots,X_n)$ and $\langle\cdot,\cdot\rangle$ is the standard inner product. We use the convention $Q_B^0(X)=|B|$. Central symmetry of $B$ implies the vanishing of all odd moments; symmetry of $Q_B^k$ in the variables $X_1,\dots,X_n$ follows when $B$ is invariant under coordinate permutations, as are Lee balls and Lee spheres. After expressing $ Q_{B^n(r)}^k(X) $ in terms of power sum symmetric functions $ s_{2i} $, writing $ \mathcal{P}(k) $ for the set of integer partitions of $ k $ and $ \lambda=(\lambda_1,\dots,\lambda_\ell) $ for a partition of $ k $, the $ Q $-polynomial has the following form,
	\begin{align}\label{eq:P1-Qp_in_powersum}
		Q_{B^n(r)}^k(X)=\sum_{\lambda\in\mathcal{P}(k)}c_\lambda \prod_{i=1}^{\ell}s_{2\lambda_i},
	\end{align}
	where $c_\lambda$ is an integer. 
	Qureshi determined $c_{(k)}$ (the coefficient of $s_{2k}$) in \eqref{eq:P1-Qp_in_powersum}
		\begin{align*}
			c_{(k)}=2\sum_{i=0}^{r}(r-i)^{2k-1}|B^n(i)|
		\end{align*}
	and proved the nonexistence of lattice tilings by Lee balls with fixed radius $ r $ in $\Z^n$ for certain infinite sets of dimensions $n$ in \cite{qureshi_non-existence_2020}. This framework covers both Kim's approach and the Zhang--Ge approach in \cite{kim_nonexistence_2017,zhang_perfect_2017}.

	For small values of $r$ and arbitrary $k$ and $n$, the coefficients of every term in \eqref{eq:P1-Qp_in_powersum} are known.  However, for large values of $r$, an explicit formula for these coefficients remains unknown.

	Besides the symmetric polynomial approach, group-ring equations also provide an effective tool for investigating the existence of lattice tilings in $\Z^n$. Moreover, they have been widely used in the study of representation theory, difference sets and related topics; we refer to \cite{beth_design_1999,lander_symmetric_1983,pott_finite_1995} for details on this topic.
	In \cite{zhang_nonexistence_2019}, necessary and sufficient conditions for lattice tilings by Lee balls of radii $2$ and $3$ were formulated in terms of group-ring equations as follows.
	\begin{lem}[\cite{zhang_nonexistence_2019}]\label{lem:P5-groupring_r2}
		Let $n\ge2$. There exists a lattice tiling of $\Z^n$ by $B^n(2)$ if and only if there exists a finite abelian group $G$ of order $2n^2+2n+1$ and a subset $T$ of size $2n+1$ viewed as an element in $\Z[G]$ satisfying
		\begin{itemize}
			\item the identity element $e$ belongs to $T$,
			\item $T=T^{(-1)}$,
			\item $T^2=2G-T^{(2)}+2n$.
		\end{itemize}
	\end{lem}
	Furthermore, necessary and sufficient conditions in terms of group-ring equations for lattice tilings of $\Z^n$ by $B^n(3)$ were also proved in \cite{zhang_nonexistence_2019}. 
	Based on them, Zhang and the second author derived necessary conditions for the existence of lattice tilings by Lee balls of radii 2 and 3 in \cite{zhang_nonexistence_2019}. In particular, the following corollary can be proved.
	\begin{cor}[\cite{zhang_nonexistence_2019}]\label{thm:P5-zhangandzhou}
		Assume that $n\equiv1,5\pmod{7}$. If $24n+1$ is not a square or $84\nmid(24n+1)^2\pm6\sqrt{24n+1}+29$, then there is no lattice tiling by $B^n(3)$.
	\end{cor}

	By computing directly with the equations in Lemma \ref{lem:P5-groupring_r2} modulo $3$ and $5$, Leung and the second author proved the Golomb--Welch conjecture in the lattice-tiling case for radius $r=2$.
	\begin{thm}[\cite{leung_no_2020}]
		For any integer $n\geq 3$, there is no lattice tiling of $\Z^n$ by Lee ball $B^n(2)$.
	\end{thm}
	This approach can also be applied to other lattice-tiling problems; see \cite{zhang_linear_2024,zhang_limited_magnitude_2023}.
	These results show that group-ring equations are effective tools for proving nonexistence. They motivate the following question.
	\begin{question}\label{que:group.ring}
		What are the necessary and sufficient conditions in terms of group-ring equations for the existence of a lattice tiling by $B^n(r)$ for arbitrary radius $r$, analogous to Lemma \ref{lem:P5-groupring_r2}?
	\end{question}

	\subsection{Our Contribution}

	In this paper, we focus on the lattice cases of the Golomb--Welch conjecture. We work primarily within the two frameworks described above, based on symmetric polynomials and group-ring equations.

	Our first main contribution is an explicit formula for every coefficient $c_\lambda$ in \eqref{eq:P1-Qp_in_powersum}, proved using enumerative-combinatorial and $q$-series methods.
	\begin{thm}\label{thm:P3-Q-poly_explicit_formula}
		Let $n\ge1$, $r\ge0$, and $k\ge1$. For every partition $\lambda=(\lambda_1,\lambda_2,\dots,\lambda_\ell)\in\mathcal P(k)$, define
		\begin{align}\label{eq:P3-cor_Q-poly_explicit_formula}
			c_\lambda
			:={}&\frac{(2k)!}{\prod_{i=1}^{\ell}(2\lambda_i)!}
			\cdot\frac{2^\ell}{\prod_{j=1}^{k}m_\lambda(j)!}
			\sum_{\substack{r_1+\cdots+r_{\ell+1}=r\\ r_1,\dots,r_\ell\ge1,\ r_{\ell+1}\ge0}}
			|B^n(r_{\ell+1})|\prod_{i=1}^{\ell}r_i^{2\lambda_i-1}.
		\end{align}
		Here $m_\lambda(j)$ denotes the multiplicity of the part $j$ in $\lambda$. Then
		\[
		Q_{B^n(r)}^k(X)=\sum_{\lambda\in\mathcal P(k)}c_\lambda s_{2\lambda}.
		\]
	\end{thm}

	Our second main contribution is a necessary and sufficient group-ring characterization of lattice tilings by $B^n(r)$ for arbitrary radius $r$, answering Question~\ref{que:group.ring}.
	\begin{thm}\label{thm:P5-groupring}
		Let $n,r\ge1$. There exists a lattice tiling of $\Z^n$ by $B^n(r)$ if and only if there exist a finite abelian group $G$ of order $|B^n(r)|$ and a subset $T\subseteq G$ of size $2n+1$, regarded as an element of $\Z[G]$, such that
		\begin{itemize}
			\item the identity element $e$ belongs to $T$;
			\item $T=T^{(-1)}$;
			\item
			\[
			G=\sum_{t=0}^{\lfloor r/2\rfloor}(-1)^t\binom{n}{t}
			\left(\sum_{\lambda\in\mathcal P(r-2t)}
			\frac{T^{(\lambda)}}{\prod_{i=1}^{\ell(\lambda)}\lambda_i
			\prod_{j=1}^{r-2t}m_\lambda(j)!}\right).
			\]
		\end{itemize}
		For $\lambda=(\lambda_1,\dots,\lambda_\ell)$, we set
		$T^{(\lambda)}=T^{(\lambda_1)}\cdots T^{(\lambda_\ell)}$, and for the empty partition we set $T^{(\varnothing)}=e$.
	\end{thm}

	To apply Theorem~\ref{thm:P3-Q-poly_explicit_formula} efficiently to specified values of $n$ and $r$, we prove the following ideal-theoretic nonexistence criterion.
	\begin{thm}\label{thm:P4-ideal}
		Let \(p>2\) be a prime. Assume \(B^n(r)\subseteq\Z^n\) is a Lee ball such that $|B^n(r)|=pm$ with \(p\nmid m\). In the multivariate polynomial ring $\F_p[Y_1,\dots,Y_{\frac{p-1}{2}}]$, if the polynomial $\sum_{\lambda\in\mathcal{P}(\frac{p-1}{2})}c_\lambda Y_{\lambda}$ lies within the ideal generated by the polynomials $\sum_{\lambda\in\mathcal{P}(k)}c_\lambda Y_{\lambda}$ for $1\le k<\frac{p-1}{2}$, i.e.,
		\[
		\sum_{\lambda\in\mathcal{P}(\frac{p-1}{2})}c_\lambda Y_{\lambda}\in\left( \sum_{\lambda\in\mathcal{P}(1)}c_\lambda Y_{\lambda}, \sum_{\lambda\in\mathcal{P}(2)}c_\lambda Y_{\lambda}, \dots, \sum_{\lambda\in\mathcal{P}(\frac{p-1}{2}-1)}c_\lambda Y_{\lambda} \right),
		\]
		where $Y_\lambda=Y_{\lambda_1}Y_{\lambda_2}\cdots Y_{\lambda_\ell}$ and $\left(f_1,\dots,f_k\right)$ denotes the ideal generated by the polynomials $f_1,\dots,f_k$, then no lattice tiling of $\Z^n$ by the Lee ball $B^n(r)$ exists.
	\end{thm}

	A consequence of Theorem~\ref{thm:P4-ideal} is that, for fixed $r$ and $p$, Theorem~\ref{thm:P3-Q-poly_explicit_formula} and symbolic computation rule out lattice tilings of $\Z^n$ by $B^n(r)$ for a fixed proportion of positive integer dimensions $n$.

	More precisely, using the computations described in Section~\ref{sec:nonexi_inf_n}, which are exhaustive for prime $p<40$ and partial for $40<p<50$, and considering all radii \(3\le r\le50\), we obtain Table~\ref{ta:proportion}. In this table, the column labeled ``\(r\)'' gives the radius of the Lee ball, and the column labeled ``proportion'' gives, for each such radius, the proportion of positive integers \(n \in \Z^+\) for which Theorem~\ref{thm:P4-ideal} can prove the nonexistence of a lattice tiling by Lee balls. Notably, Table \ref{ta:proportion} shows that, for all radii \(3 \le r \le 50\), the set of dimensions \(n\) for which Theorem~\ref{thm:P4-ideal} proves the nonexistence of such tilings has positive density in \(\Z^+\).
	\begin{table}[!ht]
		\centering
		\caption{Proportions of Dimensions $n$ for Which Theorem~\ref{thm:P4-ideal} Proves Nonexistence, for Each Radius $r$}
		\resizebox{\textwidth}{!}{%
		\begin{tabular}{cc|cc|cc|cc|cc|cc}
			\hline
			$r$ & proportion & $r$ & proportion & $r$ & proportion & $r$ & proportion & $r$ & proportion & $r$ & proportion \\ \hline
			3 & 10.49\% & 11 & 24.50\% & 19 & 20.43\% & 27 & 44.72\% & 35 & 42.83\% & 43 & 47.09\% \\
			4 & 14.81\% & 12 & 50.60\% & 20 & 8.59\% & 28 & 32.07\% & 36 & 51.76\% & 44 & 46.86\% \\
			5 & 10.37\% & 13 & 28.67\% & 21 & 13.91\% & 29 & 35.43\% & 37 & 52.79\% & 45 & 42.35\% \\
			6 & 13.79\% & 14 & 36.91\% & 22 & 34.91\% & 30 & 34.16\% & 38 & 58.39\% & 46 & 29.63\% \\
			7 & 22.97\% & 15 & 25.83\% & 23 & 18.37\% & 31 & 34.38\% & 39 & 45.82\% & 47 & 39.82\% \\
			8 & 10.16\% & 16 & 33.04\% & 24 & 22.97\% & 32 & 50.63\% & 40 & 35.41\% & 48 & 30.36\% \\
			9 & 21.33\% & 17 & 28.93\% & 25 & 27.07\% & 33 & 37.48\% & 41 & 41.78\% & 49 & 39.75\% \\
			10 & 33.74\% & 18 & 9.59\% & 26 & 19.11\% & 34 & 37.88\% & 42 & 62.66\% & 50 & 47.27\% \\ \hline
		\end{tabular}}
	\label{ta:proportion}
	\end{table}

	The final part of our contribution is a series of necessary conditions for the existence of lattice tilings by Lee balls, obtained via symmetric functions and Fourier analysis over finite fields. We then propose an algorithm that integrates these conditions to determine the nonexistence of such a tiling for any fixed dimension $n$ and radius $ r $.
	We apply it to all tabulated pairs $(n,r)$ with $3\le r\le12$ and $|B^n(r)|<2^{30}$. There are $1{,}358$ such pairs, and our algorithm with some extra analysis for special parameters rules out $1{,}356$ of them, leaving only two open cases: $(35,3)$ and $(805,3)$; see Algorithm~\ref{alg:nonexistence-tester} and Theorem~\ref{thm:main_computational}.

	Through the correspondence between lattice tilings by Lee balls and the Abelian Cayley graphs meeting the Moore-like bound, our achievements for dimension $n$ and radius $r$ also improve the corresponding extremal problem for Abelian Cayley graphs with degree $2n$ and diameter $r$. For instance, Theorem~\ref{thm:P5-groupring} gives a group-ring characterization, valid for arbitrary radius, of the equality cases in the Abelian Cayley--Moore bound. 

	\subsection{Organization of the paper} Section~\ref{sec:pre} introduces the notation and preliminaries. Section~\ref{sec:Q_poly} proves Theorem~\ref{thm:P3-Q-poly_explicit_formula}, giving the Q-polynomial expansion for Lee balls and proves the group-ring characterization in Theorem~\ref{thm:P5-groupring}. Section~\ref{sec:nonexi_inf_n} develops ideal-theoretic necessary conditions and proves positive-density nonexistence results for every radius $3\le r\le50$. Section~\ref{sec:algorithm} derives further conditions and combines them in a nonexistence-testing algorithm. For $(n,r)$ with $|B^n(r)|<2^{30}$, the combined method rules out $1{,}356$ of $1{,}358$ pairs.

	\section{Notations and Preliminaries}\label{sec:pre}
	\subsection{Lattice Tilings by Lee Balls}

	In \cite{stanton_note_1970}, Stanton and Cowan showed that the Lee ball $B^n(r)$ has size
	\begin{align}\label{eq:P2-size_ball}
		|B^n(r)|=\sum_{i=0}^{\min\{n,r\}}2^i\binom{n}{i}\binom{r}{i}.
	\end{align}
	Moreover, in \cite{post_nonexistence_1975}, Post determined its generating function
	\begin{align}\label{eq:P2-series_lee_ball}
		\sum_{r=0}^\infty|B^n(r)|q^r=\frac{(1+q)^n}{(1-q)^{n+1}}.
	\end{align}
	Define the \emph{Lee sphere} (centered at $\mathbf{0}$) by
	$$S^n(r):=\{\mathbf{x}\in\Z^n:d_L(\mathbf{x},\mathbf{0})=r\}.$$
	With the convention $B^n(-1)=\varnothing$, we have
	\begin{align}\label{eq:P2-series_lee_sphere}
		\sum_{r=0}^\infty|S^n(r)|q^r=&\sum_{r=0}^\infty\left(|B^n(r)|-|B^n(r-1)|\right)q^r\notag\\
		&=\sum_{r=0}^\infty|B^n(r)|q^r-q\sum_{r=0}^\infty|B^n(r)|q^r=\left(\frac{1+q}{1-q}\right)^n.
	\end{align}

	The following group-homomorphism criterion for lattice tilings of $\Z^n$ was proved in \cite{horak_diameter_2012}.

	\begin{lem}[\cite{horak_diameter_2012}]
		Let $ S $ be a finite subset of \( \mathbb{Z}^n \) with \( |S| = m \). Then there exists a lattice tiling of \( \mathbb{Z}^n \) by \( S \) if and only if there exists an abelian group \( G \) of order \( m \) and a homomorphism \( \varphi: \mathbb{Z}^n \to G \) that is bijective on \( S \).
	\end{lem}
	As an immediate consequence, we obtain the following corollary for Lee balls.
	\begin{cor}[\cite{horak_diameter_2012}]
		\label{cor:P2-group_latticetiling}
		There exists a lattice tiling of \( \mathbb{Z}^n \) by the Lee ball of radius \( r \) if and only if there exists an abelian group \( G \) of size \( |B^n(r)| \) and a homomorphism \( \varphi: \mathbb{Z}^n \to G \) that is bijective on the ball \( B^n(r) \).
	\end{cor}

	\subsection{Integer Partitions, Symmetric Functions and Complete Bell Polynomials}
	A \emph{partition of a positive integer $n$} is a finite sequence of weakly decreasing positive integers $\lambda=(\lambda_1,\lambda_2,\dots,\lambda_\ell)$ such that $n=\lambda_1+\lambda_2+\cdots+\lambda_\ell$. The size of the partition is denoted by $|\lambda|(=n)$. Each $\lambda_i$ is called a \emph{part} of the partition $\lambda$. The number of parts is called the \emph{length} of the partition, denoted by $\ell(\lambda)$. The number of times a part $\lambda_i$ appears in the partition is called the \emph{multiplicity} of $\lambda_i$, denoted by $m_\lambda(\lambda_i)$. Sometimes, a partition $\lambda$ can also be represented in the following form:
	\[
	\lambda=\left(\lambda_1^{m_\lambda(\lambda_1)},\lambda_2^{m_\lambda(\lambda_2)},\dots,\lambda_k^{m_\lambda(\lambda_k)}\right),
	\]
	where $\lambda_1>\lambda_2>\cdots>\lambda_k$, and the exponent can be omitted when $m_\lambda(\lambda_i)=1$.

	We also use the standard empty-partition convention $\mathcal P(0)=\{\varnothing\}$, with $\ell(\varnothing)=0$ and every empty product equal to $1$.

	In this paper, we denote the set of all partitions of $n$ by $\mathcal{P}(n)$ and the number of partitions of $n$ by $p(n)$. For example,
	\[
	\mathcal{P}(5) = \{(5), (4,1), (3,2), (3,1,1), (2,2,1), (2,1,1,1), (1,1,1,1,1)\}
	\]
	and $p(5) = |\mathcal{P}(5)| = 7$. For a positive integer $k$ and a partition $\lambda=(\lambda_1,\ldots,\lambda_\ell)$, we write $k\lambda=(k\lambda_1,\ldots,k\lambda_\ell)$.

	A \emph{symmetric polynomial} is invariant under permutations of its variables. We use two standard families: elementary symmetric polynomials and power-sum symmetric polynomials.
	\begin{defn}[Elementary Symmetric Polynomials]
		Let $X=(X_1,\dots,X_n)$. For an integer \(k\) with $1\le k\le n$, the \(k\)-th elementary symmetric polynomial \( e_k(X) \) is defined to be
		$$e_k(X):=\sum_{1\le i_1<\cdots<i_k\le n}X_{i_1}\cdots X_{i_k}.$$
	\end{defn}

	\begin{defn}[Power Sum Symmetric Polynomials]
		Let $X=(X_1,\dots,X_n)$. For a positive integer \(k\), the \(k\)-th power-sum symmetric polynomial \( s_k(X) \) is defined to be
		\[
		s_k(X) = \sum_{i=1}^n X_i^k.
		\]
	\end{defn}

	For a partition $\lambda=(\lambda_1,\dots,\lambda_\ell)$, we use the standard notation
	$$e_\lambda:=e_{\lambda_1}\cdots e_{\lambda_\ell}\text{ and }s_\lambda:=s_{\lambda_1}\cdots s_{\lambda_\ell}.$$
	
	The \emph{complete Bell polynomials} \cite{riordan_introduction_2014} \(B_n(a_1,\dots,a_n)\) are defined through the exponential generating function
	\begin{align}\label{DefnOfBellPoly}
		\exp\!\left( \sum_{i=1}^{\infty} a_i \frac{x^i}{i!} \right) = \sum_{n=0}^{\infty} B_n(a_1,\dots,a_n) \frac{x^n}{n!},
	\end{align}
	with \(B_0 = 1\).  Equivalently, they admit the explicit combinatorial form
	\begin{align}\label{ExpandBellPoly}
	B_n(a_1,\dots,a_n) = \sum_{\lambda \in \mathcal{P}(n)} \frac{n!}{\prod_{i} \lambda_i!} \frac{1}{\prod_{j} m_{\lambda}(j)!} \prod_{i=1}^{\ell(\lambda)} a_{\lambda_i},
	\end{align}
	where the sum runs over all integer partitions \(\lambda\) of \(n\), \(\lambda_i\) are the parts of \(\lambda\), \(\ell(\lambda)\) is the number of parts, and \(m_{\lambda}(j)\) is the multiplicity of \(j\) in \(\lambda\).  This identity is fundamental in combinatorics and probability, encoding the relation between moments and cumulants.

	\subsection{Group-Ring Equations on Abelian Groups and Their Fourier Transforms}

	Let $G$ be a group and $R$ a ring. The group ring $R[G]$ consists of the finitely supported functions $G\to R$ and is the free $R$-module with basis $\{g:g\in G\}$. When $G$ is finite, every function $G\to R$ may be identified with a group-ring element. 
	
	For a finite multiset $A$ on $G$, we identify $A$ with the group-ring element
	$\sum_{g\in G}a_g g\in\Z[G]$, where $a_g$ is the multiplicity of $g$ in $A$. Given any \( A = \sum_{g \in G} a_g g \in \mathbb{Z}[G] \), we define
	$A^{(t)} = \sum_{g \in G} a_g g^t.$
	Addition and multiplication in group rings are defined by:
	\[
	\sum_{g\in G}a_gg+\sum_{g\in G}b_gg=\sum_{g\in G}(a_g+b_g)g,
	\]
	and
	\[
	\sum_{g\in G}a_gg \sum_{g\in G}b_gg=\sum_{g\in G}\left(\sum_{h\in G}a_hb_{h^{-1}g}\right)g.
	\]

	For a field $\F$, write $\F^*$ for its multiplicative group. A \emph{character} of a group $G$ is a homomorphism $\chi:G\rightarrow\F^*$, i.e., an element of the $\Z$-module $\operatorname{Hom}_\Z(G,\F^*)$. When $G$ is finite, a character may also be identified with an element of $\F[G]$.

	The characters of $G$ form the \emph{character group}, denoted by $\hat{G}$. The identity element of $\hat{G}$ is referred to as the \emph{principal character} of $G$, denoted by $1^G$ (sometimes also denoted by $\chi_0$).

	In this paper, we only have to handle the characters of finite abelian groups $G$. Given a function $f:G\rightarrow\F$, its \emph{Fourier transform} is the function $\hat{f}:\hat{G}\rightarrow\F$ defined by
	\begin{equation}\label{eq:def_fourier}
		\hat{f}(\chi):=\sum_{g\in G}f(g)\chi(g^{-1}),\ \chi\in\hat{G}.
	\end{equation}
	The function $f$ can be recovered from its Fourier transform by the \emph{inversion formula}
	$$f(g)=|G|^{-1}\sum_{\chi\in\hat{G}}\hat{f}(\chi)\chi(g),\ g\in G.$$
	See \cite[Chapter 6]{beth_design_1999} and \cite{ceccerini_discrete_2018} for more results on Fourier transforms.

	The support of $ f $ is defined by
	$$\mathrm{supp}(f)=\{x\in G:f(x)\neq0\}.$$
	In \cite{feng_shift_2019}, Feng, Hollmann and Xiang proved the following uncertainty principle on the support of $f$ and $\hat{f}$.
	\begin{thm}[\cite{feng_shift_2019}]\label{thm:P6-uncertainty_principle}
		Let $ (G,+) $ be a finite abelian group and let $ \F $ be a field of characteristic zero or characteristic $ p $ not dividing $ |G| $. Let $ \mathbb{E} $ denote an extension of $ \F $ containing a primitive $ N $-th root of unity, where $ N =\exp(G) $ is the exponent of $G$, and let $ \hat{G} $ denote the group of $ \mathbb{E} $-valued characters of $ G $. Then for any nonzero function $ f : G \rightarrow \F $ with Fourier transform $ \hat{f} $, we have that
		\[
		|\mathrm{supp}(f)||\mathrm{supp}(\hat{f})|\geq|G|+|\mathrm{supp}(\hat{f})|-|H(\mathrm{supp}(\hat{f}))|;
		\]
		and dually,
		\[
		|\mathrm{supp}(f)||\mathrm{supp}(\hat{f})|\geq|G|+|\mathrm{supp}({f})|-|H(\mathrm{supp}({f}))|,
		\]
		where $H(\mathrm{supp}({f}))$ is the stabilizer of $\mathrm{supp}(f)$ in $ G $.
	\end{thm}
	By Theorem \ref{thm:P6-uncertainty_principle}, one directly obtains
	\begin{align}\label{eq:P6-UP}
		|\mathrm{supp}(f)||\mathrm{supp}(\hat{f})|\geq|G|.
	\end{align}

	Later we concentrate on two classes of finite abelian groups and their subgroups: the additive group of a finite field $(\F_q,+)$ and the multiplicative group of its nonzero elements $(\F^*_q,\times)$.

	By the fundamental theorem of finite abelian groups, the following lemma can be directly proved.
	\begin{lem}\label{lem:P2-m_to_1}
		Let $G$ be a finite abelian group of order $pm$, where $p$ is prime. Then there exists a surjective homomorphism
		$\Phi:G\rightarrow(\F_p,+)$. In particular, $\Phi$ is an $m$-to-$1$ map.
	\end{lem}

	Now consider the multiplicative group $\F_p^*$ which is of order $p-1$. For any subgroup $H\le\F_p^*$, since $H$ is cyclic (as a subgroup of a cyclic group), we can explicitly write all its characters over $\F_p$ as
	\[
	\chi_k: H\to\F_p^*,\ \chi_k(h)=h^k\ (k=0,1,\dots,|H|-1),
	\]
	where $h^0=1$ (the principal character). The set $\hat H=\{\chi_0,\dots,\chi_{|H|-1}\}$ is the character group of $H$, and $|\hat H|=|H|$.

	By \eqref{eq:def_fourier}, for a function $f:H\to\F_p$, its Fourier transform $\hat f:\hat H\to\F_p$ is given by
	\[
	\hat f(\chi_k):=\sum_{h\in H}f(h)\,\chi_k(h^{-1})=\sum_{h\in H}f(h)\,h^{-k},\quad k=0,\dots,|H|-1.
	\]
	The following theorem summarizes two fundamental properties of these characters which follow directly from orthogonality and the Fourier inversion formula.
	\begin{thm}
		\label{thm:P6-Fourier_reverse}
		Let $\chi_s:H\rightarrow\F_p^*$ be a character of $H$, where $H$ is a subgroup of $(\F_p^*,\times)$.
		\begin{enumerate}[label=(\alph*)]
			\item For any $a\neq1$ in $H$, $\sum_{s=0}^{|H|-1}\chi_s(a)=0$.
			\item For any function $g:H\rightarrow\F_p$,
			$$g(x)=\frac{1}{|H|}\sum_{s=0}^{|H|-1}\hat{g}(\chi_s)\chi_s(x).$$
		\end{enumerate}
	\end{thm}

	By Theorem \ref{thm:P6-Fourier_reverse}(a), for every integer $k\ge1$, we have
	\begin{equation}\label{eq:power_sum_value}
		\sum_{\alpha\in\F_p}\alpha^k=
		\begin{cases}
			-1,&\text{if }p-1\mid k;\\
			0,&\text{if }p-1\nmid k.
		\end{cases}
	\end{equation}

	\section{Explicit Formulas for Q-Polynomials and Group-Ring Equations}\label{sec:Q_poly}

	In \cite{qureshi_non-existence_2020}, Qureshi introduced a family of homogeneous polynomials.
	\begin{defn}[\cite{qureshi_non-existence_2020}]\label{def:Q-poly}
		Let $B\subseteq\Z^n$ be finite and let $k\in\Z_{\ge0}$. The \emph{Q-polynomial} of order $k$ associated with $B$ is
		\[
		Q_B^k(X):=\sum_{b\in B}\langle X,b\rangle^{2k},
		\qquad X=(X_1,\dots,X_n).
		\]
		In particular, $Q_B^0(X)=|B|$.
	\end{defn}
	Recall that $s_k=\sum_{i=1}^nX_i^k$, that $s_\lambda=s_{\lambda_1}\cdots s_{\lambda_\ell}$ for a partition $\lambda=(\lambda_1,\dots,\lambda_\ell)$, and that $k\lambda=(k\lambda_1,\dots,k\lambda_\ell)$. We will use the following observations.
	\begin{itemize}
		\item For every $k\in\Z_{\ge0}$, central symmetry gives
		$\sum_{b\in B^n(r)}\langle X,b\rangle^{2k+1}
		=\sum_{b\in S^n(r)}\langle X,b\rangle^{2k+1}=0$.
		\item Coordinate-permutation invariance implies that both $Q_{B^n(r)}^k(X)$ and $Q_{S^n(r)}^k(X)$ are symmetric polynomials.
		\item For every $k\ge1$, the coefficients $c_\lambda$ of $s_{2\lambda}$ in 
		$Q_{B^n(r)}^k(X)=\sum_{\lambda\in\mathcal P(k)}c_\lambda s_{2\lambda}$ are integers.
	\end{itemize}

	For example, Kim \cite{kim_nonexistence_2017} proved the following identity.
	\begin{align*}
	Q_{B^n(2)}^k(X)
	={}&\sum_{i=1}^{n}\sum_{\varepsilon=\pm1}
	\bigl((\varepsilon X_i)^{2k}+(2\varepsilon X_i)^{2k}\bigr)
	+\sum_{1\le i<j\le n}\sum_{\varepsilon,\delta=\pm1}
	(\varepsilon X_i+\delta X_j)^{2k}\\
	={}&(4^k+4n+2)s_{2k}
	+2\sum_{t=1}^{k-1}\binom{2k}{2t}s_{2t}s_{2(k-t)}.
	\end{align*}

	The first part of this section is devoted to the proof of Theorem \ref{thm:P3-Q-poly_explicit_formula}. To this end, we first prove Theorem \ref{thm:P3-Q-poly_sphere} and then establish the following chain of implications:
	\[
	\text{Theorem \ref{thm:P3-Q-poly_sphere}} \Rightarrow \text{Theorem \ref{thm:P3-Q-poly_ball}} \Rightarrow \text{Theorem \ref{thm:P3-Q-poly_explicit_formula}}.
	\]

	Recall that \( S^n(r) \) denotes the Lee sphere of dimension \( n \) and radius \( r \) centered at $\mathbf{0}$ and $s_\lambda$ denotes the product $s_{\lambda_1}\cdots s_{\lambda_\ell}$ of the power-sum symmetric polynomials for $\lambda = (\lambda_1,\dots,\lambda_\ell)$.

	\begin{thm}\label{thm:P3-Q-poly_sphere}
		Let \( Q_{S^n(r)}^k(X) \)  be the \( k \)-th order Q-polynomial associated with \( S^n(r) \).
		When \( k = 0 \),
		\begin{align*}
			\sum_{r=0}^\infty Q_{S^n(r)}^0(X)q^r= \left( \frac{1+q}{1-q} \right)^n.
		\end{align*}

		When \( k \ge 1 \),
		\begin{align}\label{eq:P3-Q-poly_sphere_k>0}
			&\sum_{r=0}^\infty Q_{S^n(r)}^k(X)q^r \notag\\
			&= \left( \frac{1+q}{1-q} \right)^n \sum_{\lambda \in \mathcal{P}(k)} \left[ \frac{(2k)!}{\prod_{i=1}^{\ell(\lambda)} (2\lambda_i)!} \cdot \frac{2^{{\ell(\lambda)}}}{\prod_{i=1}^k m_\lambda(i)!} \cdot \prod_{i=1}^{\ell(\lambda)} \left( \sum_{s=1}^\infty s^{2\lambda_i - 1} q^s \right) \cdot s_{2\lambda} \right].
		\end{align}

	\end{thm}

	\begin{proof}
		Consider the formal generating function
		\begin{align*}
		F_n(X;z,q)
		&:=\sum_{r\ge0}\sum_{b\in S^n(r)}q^r e^{z\langle X,b\rangle}\\
		&=\sum_{b\in\Z^n}q^{\|b\|_1}e^{z\langle X,b\rangle}=\sum_{b\in\Z^n} q^{\sum_{i=1}^{n} |b_i|} e^{z \sum_{i=1}^{n} X_i b_i}\\
		&=\sum_{b\in\Z^n}\prod_{i=1}^{n}\left( q^{|b_i|} e^{z X_i b_i} \right)= \prod_{i=1}^{n} \sum_{b_i\in\Z} \left( q^{|b_i|} e^{z X_i b_i} \right)\\
		&=\prod_{i=1}^n\left(1+\sum_{b_i\ge1}\left( q^{b_i} e^{z X_i b_i}+q^{b_i} e^{-z X_i b_i} \right)\right)\\
		&=\prod_{i=1}^n\left(1+\frac{qe^{zX_i}}{1-qe^{zX_i}}+\frac{qe^{-zX_i}}{1-qe^{-zX_i}}\right)\\
		&=\prod_{i=1}^n\left(\frac{1-q^2}{(1-qe^{zX_i})(1-qe^{-zX_i})}\right).
		\end{align*}
		For each $i$, we can obtain
		\begin{align}
		\nonumber
		\ln\left(\frac{1-q^2}{(1-qe^{zX_i})(1-qe^{-zX_i})}\right)
		&=\ln(1-q^2)+\sum_{s\ge1}\frac{q^s}{s}(e^{szX_i}+e^{-szX_i})\\
		\label{eq:ln(1-xq)}
		&=\ln(1-q^2)+2\sum_{j\ge0}\left(\sum_{s\ge1}s^{2j-1}q^s\right)\frac{(zX_i)^{2j}}{(2j)!}\\
		&=\ln\left(\frac{1+q}{1-q}\right)+2\sum_{j\ge1}\left(\sum_{s\ge1}s^{2j-1}q^s\right)\frac{(zX_i)^{2j}}{(2j)!}.\notag
		\end{align}
		Taking the sum over $ i $ gives
		\begin{equation}\label{eq:F_n.def}
			F_n(X;z,q)=\left(\frac{1+q}{1-q}\right)^n
			\exp\!\left(
			2\sum_{j\ge1}\left(\sum_{s\ge1}s^{2j-1}q^s\right)
			\frac{s_{2j}(X)z^{2j}}{(2j)!}
			\right).
		\end{equation}

		Evaluating $F_n(X;z,q)$ at $z=0$, we get the asserted formula for $k=0$. For $k\ge1$, by computing the partial derivatives of $		F_n(X;z,q)
		=\sum_{r\ge0}\sum_{b\in S^n(r)}q^r e^{z\langle X,b\rangle}$ with respect to $z$ at $z=0$, we obtain
		\[
		\sum_{r\ge0}Q_{S^n(r)}^k(X)q^r=\sum_{r\ge0} \sum_{b\in S^n(r)}\langle X,b\rangle^{2k}q^r
		=\left.\frac{\partial^{2k}}{\partial z^{2k}}F_n(X;z,q)\right|_{z=0}.
		\]
		
		By comparing \eqref{eq:F_n.def} and the exponential generating function of the complete Bell polynomials given in \eqref{DefnOfBellPoly}, we immediately find that the generating function \(F_n(X; z,q)\) can be rewritten in terms of the complete Bell polynomials as
		\[
		F_n(X; z, q) = \left( \frac{1+q}{1-q} \right)^n \sum_{m \ge 0} B_m(a_1, \dots, a_m) \frac{z^m}{m!}, 
		\]
		where
		\[
		a_{2j-1} = 0, \quad \text{and} \quad a_{2j} = 2\left( \sum_{s \ge 1} s^{2j-1} q^s \right) s_{2j}(X), \quad \text{for } j \ge 1. 
		\]

		Consequently, taking the $ 2k $-th partial derivative of $ F_n(X; z, q)$ with respect to $z$ at $ z=0 $ cancels exactly the factorial denominator $ (2k)! $, yielding
		\[
		\left. \frac{\partial^{2k}}{\partial z^{2k}} F_n(X; z, q) \right|_{z=0} = \left( \frac{1+q}{1-q} \right)^n B_{2k}(a_1, \dots, a_{2k}). 
		\]

		Applying the explicit combinatorial formula \eqref{ExpandBellPoly} for the complete Bell polynomial \(B_{2k}\), we obtain \eqref{eq:P3-Q-poly_sphere_k>0}.
	\end{proof}

	\begin{thm}\label{thm:P3-Q-poly_ball}
		Let \( Q_{B^n(r)}^k(X) \) be the \( k \)-th order Q-polynomial associated with \( B^n(r) \).
		For $k=0$,
		\begin{align*}
			\sum_{r=0}^\infty Q_{B^n(r)}^0(X)q^r=\frac{(1+q)^n}{(1-q)^{n+1}}.
		\end{align*}

		For $k \geq 1$,
		\begin{align*}
			&\sum_{r=0}^\infty Q_{B^n(r)}^k(X)q^r \notag \\
			&=\frac{(1+q)^n}{(1-q)^{n+1}}\sum_{\lambda \in \mathcal{P}(k)} \left[ \frac{(2k)!}{\prod_{i=1}^{\ell(\lambda)} (2\lambda_i)!} \cdot \frac{2^{\ell(\lambda)}}{\prod_{i=1}^k m_\lambda(i)!} \cdot \prod_{i=1}^{\ell(\lambda)} \left( \sum_{s=1}^\infty s^{2\lambda_i-1} q^s \right) \cdot s_{2\lambda} \right].
		\end{align*}
		where $\mathcal{P}(k)$ is the set of all partitions of $k$, $\lambda = (\lambda_1, \lambda_2, \dots, \lambda_\ell)$ is a partition of $k$, $m_\lambda(i)$ is the multiplicity of $i$ in $\lambda$, and $s_{2\lambda} = s_{2\lambda_1} \cdots s_{2\lambda_\ell}$.
	\end{thm}

	\begin{proof}
		Using $B^n(r)=\bigsqcup_{j=0}^r S^n(j)$, and setting $B^n(-1)=\varnothing$, we have
		\[
		Q_{S^n(r)}^k(X)=Q_{B^n(r)}^k(X)-Q_{B^n(r-1)}^k(X),
		\]
		which implies
		\begin{align*}
			\sum_{r=0}^{\infty} Q_{S^n(r)}^k(X) q^r
			&= \sum_{r=0}^{\infty} \left( Q_{B^n(r)}^k(X) - Q_{B^n(r-1)}^k(X) \right) q^r \\
			&= \sum_{r=0}^{\infty} Q_{B^n(r)}^k(X) q^r - q \sum_{r=0}^{\infty} Q_{B^n(r-1)}^k(X) q^{r-1} \\
			&= (1 - q) \sum_{r=0}^{\infty} Q_{B^n(r)}^k(X) q^r.
		\end{align*}
		Dividing the formula in Theorem~\ref{thm:P3-Q-poly_sphere} by $1-q$ gives exactly the asserted formula for Lee balls, including the case $k=0$.
	\end{proof}

	Finally, we prove Theorem~\ref{thm:P3-Q-poly_explicit_formula}, which gives an explicit expression for every coefficient $c_\lambda$ in
	\[
	Q_{B^n(r)}^k(X) = \sum_{\lambda \in \mathcal{P}(k)} c_\lambda s_{2\lambda}.
	\]

	\begin{proof}[{Proof of Theorem \ref{thm:P3-Q-poly_explicit_formula}}]
		By Theorem \ref{thm:P3-Q-poly_ball}, \( Q_{B^n(r)}^k(X) \) is the coefficient of \( q^r \) in the expansion of
		\begin{align}\label{eq:P3-proofof1.3}
			\frac{(1+q)^n}{(1-q)^{n+1}} \sum_{\lambda \in \mathcal{P}(k)} \left[ \frac{(2k)!}{\prod_{i=1}^\ell (2\lambda_i)!} \cdot \frac{2^{\ell}}{\prod_{i=1}^k m_\lambda(i)!} \cdot \prod_{i=1}^\ell \left( \sum_{s=1}^\infty s^{2\lambda_i - 1} q^s \right) \cdot s_{2\lambda} \right],
		\end{align}
		and we know that
		\[
		\frac{(1+q)^n}{(1-q)^{n+1}} = \sum_{r=0}^{\infty} |B^n(r)| q^r.
		\]
		It follows that
		\[
		Q_{B^n(r)}^k(X) = \sum_{\lambda \in \mathcal{P}(k)} \frac{(2k)!}{\prod_{i=1}^{\ell} (2\lambda_i)!} \cdot \frac{2^\ell}{\prod_{i=1}^{k} m_\lambda(i)!} \sum_{\substack{r_1+\cdots+r_{\ell+1}=r\\ r_1,\dots,r_\ell\ge1,\ r_{\ell+1}\ge0}}
			\left(|B^n(r_{\ell+1})|\prod_{i=1}^{\ell}r_i^{2\lambda_i-1}\right)s_{2\lambda}.
		\]
		{This proves \eqref{eq:P3-cor_Q-poly_explicit_formula}.}
	\end{proof}

	Next, we use enumerative-combinatorial tools to prove the necessary and sufficient group-ring condition in Theorem~\ref{thm:P5-groupring}.

	\begin{proof}[Proof of Theorem \ref{thm:P5-groupring}]
		Suppose first that a lattice tiling exists. By Corollary~\ref{cor:P2-group_latticetiling}, there are a finite abelian group $G$ of order $|B^n(r)|$ and a homomorphism
		$\varphi:\Z^n\to G$ whose restriction to $B^n(r)$ is bijective. Put
		$x_i=\varphi(\epsilon_i)$ and
		\[
		T=e+\sum_{i=1}^n(x_i+x_i^{-1}).
		\]
		For $s\ge0$, let
		\[
		L_s:=\sum_{|a_1|+\cdots+|a_n|\le s}x_1^{a_1}\cdots x_n^{a_n}\in\Z[G].
		\]
		A direct calculation in the formal power-series ring $\Z[G][[q]]$ gives
		\begin{align*}
		\sum_{s\ge0}L_sq^s
		&=\frac1{1-q}\sum_{s\ge0}\left(\sum_{|a_1|+\cdots+|a_n|=s}x_1^{a_1}\cdots x_n^{a_n}\right)q^s\\
		&=\frac1{1-q}\prod_{i=1}^n
		\left(1+\sum_{a\ge1}(x_i^a+x_i^{-a})q^a\right)\\
		&=\frac1{1-q}\prod_{i=1}^n\frac{1-q^2}{(1-x_iq)(1-x_i^{-1}q)}
		\end{align*}
		By a method similar to that of \eqref{eq:ln(1-xq)},
		\begin{align*}
		\sum_{s\ge0}L_sq^s
		&=\frac1{1-q}\prod_{i=1}^n\frac{1-q^2}{(1-x_iq)(1-x_i^{-1}q)}\\
		&=(1-q^2)^n\exp\!\left(\sum_{i=1}^{n}\sum_{j\ge1}\frac{q^j}{j}(1+x_i^j+x_i^{-j})\right)\\
		&=(1-q^2)^n
		\exp\!\left(\sum_{j\ge1}\frac{T^{(j)}q^j}{j}\right).
		\end{align*}
		
		Using the exponential formula for partitions in \eqref{DefnOfBellPoly} and \eqref{ExpandBellPoly},
		\[
		\exp\!\left(\sum_{j\ge1}\frac{T^{(j)}q^j}{j}\right)
		=\sum_{s\ge0}B_s\left(T^{(1)},T^{(2)},\dots,T^{(s)}\right)q^s
		=\sum_{s\ge0}\left(
		\sum_{\lambda\in\mathcal P(s)}
		\frac{T^{(\lambda)}}{\prod_{i=1}^{\ell(\lambda)}\lambda_i
		\prod_{j=1}^{s}m_\lambda(j)!}
		\right)q^s.
		\]
		Taking the coefficient of $q^r$ yields the displayed group-ring identity because $L_r=G$.

		Conversely, suppose that $G$ and $T$ satisfy the three stated conditions. The integer $|B^n(r)|$ is odd, so $G$ has no nonidentity involution. Hence, after choosing one representative from each inverse pair, the inverse-closed set $T$ can be written as
		\[
		T=\{e,x_1,x_1^{-1},\dots,x_n,x_n^{-1}\}.
		\]
		Define $\varphi:\Z^n\to G$ by $\varphi(\epsilon_i)=x_i$. The same formal identity shows that the right-hand side in the theorem is precisely
		$\sum_{b\in B^n(r)}\varphi(b)$. Since it equals $G$ in $\Z[G]$, every element of $G$ occurs with coefficient one. Thus $\varphi|_{B^n(r)}$ is bijective, and Corollary~\ref{cor:P2-group_latticetiling} gives the desired lattice tiling.
	\end{proof}

	As a scalar specialization of the generating-function identity in the proof of Theorem \ref{thm:P5-groupring}, replace every $T^{(j)}$ by an integer $k$. Then
		\begin{align}\label{eq:P5-generate_groupring}
			\sum_{\lambda\in\mathcal{P}(r)}\frac{k^{\ell(\lambda)}}{\prod_{i=1}^{\ell(\lambda)}\lambda_i\prod_{j=1}^{r}m_\lambda(j)!}
			=&[q^r]\prod_{j=1}^{\infty}\left(\sum_{m_\lambda(j)=0}^{\infty}\frac{1}{m_\lambda(j)!}\left(\frac{ k\cdot q^j}{j}\right)^{m_\lambda(j)}\right)  \notag \\
			=&[q^r]\prod_{j=1}^{\infty}\exp\left(\frac{k\cdot q^j}{j}\right)=[q^r]\exp\left(k\cdot\sum_{j=1}^{\infty}\frac{q^j}{j}\right)\\\notag
			=&[q^r]\frac{1}{(1-q)^k}=\binom{k+r-1}{r},
		\end{align}
		where $ [q^r]f(q) $ denotes the coefficient of $ q^r $ in the formal power series expansion of $ f(q) $.  We will use \eqref{eq:P5-generate_groupring} in Subsection \ref{sub:IntractableCases}.

	\section{Nonexistence Results for Infinitely Many Dimensions $n$}\label{sec:nonexi_inf_n}
	The additional information about Q-polynomials obtained in Section~\ref{sec:Q_poly}, which goes beyond the description in \cite{qureshi_non-existence_2020}, allows us to derive stronger nonexistence results. Recall that $ c_\lambda $ is defined in \eqref{eq:P1-Qp_in_powersum} and determined in \eqref{eq:P3-cor_Q-poly_explicit_formula}. The main result of this section is the ideal-theoretic criterion in Theorem~\ref{thm:P4-ideal}, which gives congruence constraints on $n$ that rule out lattice tilings by Lee balls.

	We will prove Lemma \ref{thm:P4-equations_system_X}, and then establish the following chain of implications:
		\begin{align}\label{eq:P5-}
			\text{Lemma \ref{thm:P4-equations_system_X}}\Rightarrow\text{Theorem \ref{thm:P4-equations_system_S}}\Rightarrow\text{Theorem \ref{thm:P4-ideal}}\Rightarrow\text{Theorem \ref{thm:P4-p-condition}}.
	\end{align}

	Note that the following criterion obtained by Qureshi in \cite{qureshi_non-existence_2020}, called the ``$p$-condition'', is a direct consequence of Theorem \ref{thm:P4-ideal}.

	\begin{thm}[\cite{qureshi_non-existence_2020}]\label{thm:P4-p-condition}
		Let \(p>2\) be a prime. Assume that the Lee ball \(B^n(r)\) satisfies \(|B^n(r)|=mp\) with \(p\nmid m\) and the following congruences are satisfied:
		\begin{align*}
			\begin{cases}
				c_{(k)}\not\equiv0\pmod{p} \text{ for }1\le k<\frac{p-1}{2};\\
				c_{(\frac{p-1}{2})}\equiv0\pmod{p}.
			\end{cases}
		\end{align*}
		Then there is no lattice tiling of \(\Z^n\) by \(B^n(r)\).
	\end{thm}

	Corollary~\ref{cor:P2-group_latticetiling} reduces the tiling problem to the existence of an $\F_p$-rational point on a certain algebraic variety, as made precise in the following lemma.
	\begin{lem}\label{thm:P4-equations_system_X}
		{Let $p>2$ be a prime. Assume that the Lee ball \(B^n(r)\) satisfies \(|B^n(r)|=mp\) and there is no $X=(X_1,\dots,X_n)\in\F_p^n$ such that the following system of equations holds:}
		\begin{align}\label{eq:P4-thm_equations_system_X}
			\begin{cases}
				Q_{B^n(r)}^k(X)\equiv0\pmod{p},\text{ for }1\le k<\frac{p-1}{2};\\
				Q_{B^n(r)}^{\frac{p-1}{2}}(X)\equiv-m\pmod{p}.
			\end{cases}
		\end{align}
		Then there is no lattice tiling of \(\Z^n\) by \(B^n(r)\).
	\end{lem}

	\begin{proof}[Proof of Lemma \ref{thm:P4-equations_system_X}]
		We prove the contrapositive. Suppose that there exists a lattice tiling of \(\Z^n\) by \(B^n(r)\). By Corollary \ref{cor:P2-group_latticetiling}, there is an abelian additive group \(G\) and a homomorphism \(\varphi_+:\Z^n\rightarrow G\) such that the restriction \(\varphi_+|_{B^n(r)}:B^n(r)\rightarrow G\) is a bijection.
		Because $p\mid|B^n(r)|$, Lemma~\ref{lem:P2-m_to_1} gives a surjective homomorphism $\Phi:G\twoheadrightarrow(\F_p,+)$. Set $x_i=(\Phi\circ\varphi_+)(e_i)$, where $e_1,\dots,e_n$ are the standard basis vectors of $\Z^n$. Then
		\[
		(\Phi\circ\varphi_+)(B^n(r))
		=\{*~a_1x_1+\cdots+a_nx_n:
		a_i\in\Z,\ |a_1|+\cdots+|a_n|\le r~*\}
		=m\F_p.
		\]

		By \eqref{eq:power_sum_value}, we obtain
		\[
		\sum_{\alpha\in\Phi\circ\varphi_+(B^n(r))}\alpha^{2k}=
		\begin{cases}
			0,&\text{if }p-1\nmid 2k;\\
			-m,&\text{if }p-1\mid 2k.
		\end{cases}
		\]

		Meanwhile, notice that \(\Phi\circ\varphi_+(B^n(r))=\{*~\left\langle x,b\right\rangle:b\in B^n(r)~*\}\), where \(x=(x_1,\dots,x_n)\in\F_p^n\), i.e.,
		\[
		\sum_{b\in B^n(r)}\left\langle x,b\right\rangle^{2k}=
		\begin{cases}
			0,&\text{if }p-1\nmid 2k;\\
			-m,&\text{if }p-1\mid 2k.
		\end{cases}
		\]

		Since \(Q_{B^n(r)}^k(X)=\sum_{b\in B^n(r)}\left\langle X,b\right\rangle^{2k}\), \(x=(x_1,\dots,x_n)\) is a solution to the system \eqref{eq:P4-thm_equations_system_X}.
	\end{proof}

	\begin{thm}\label{thm:P4-equations_system_S}
		Let \(p>2\) be a prime. Assume that the Lee ball \(B^n(r)\) satisfies $|B^n(r)|=mp$ and there do not exist \(Y_1,\dots,Y_{\frac{p-1}{2}}\in\F_p\) such that the following system of equations holds:
		\begin{align}\label{eq:P4-cor_equations_system_S}
			\begin{cases}
				&\sum_{\lambda\in\mathcal{P}(k)}c_\lambda Y_{\lambda}\equiv0\pmod{p},\text{ for }1\le k<\frac{p-1}{2};\\
				&\sum_{\lambda\in\mathcal{P}(\frac{p-1}{2})}c_\lambda Y_{\lambda}\equiv-m\pmod{p},
			\end{cases}
		\end{align}
		where \(Y_{\lambda}=Y_{\lambda_1}\cdots Y_{\lambda_\ell}\) for \(\lambda=(\lambda_1,\dots,\lambda_\ell)\) and $c_\lambda$ is defined in \eqref{eq:P1-Qp_in_powersum}. Then there is no lattice tiling of \(\Z^n\) by \(B^n(r)\).
	\end{thm}
	\begin{proof}[Proof of the implication chain \eqref{eq:P5-}]
		For the first implication, we prove the contrapositive. Suppose that a lattice tiling exists. Lemma~\ref{thm:P4-equations_system_X} gives a solution $x=(x_1,\dots,x_n)\in\F_p^n$ of \eqref{eq:P4-thm_equations_system_X}. Since
		\[
		Q_{B^n(r)}^k(X)=\sum_{\lambda\in\mathcal P(k)}c_\lambda s_{2\lambda}(X),
		\]
		the tuple
		\[
		(Y_1,\dots,Y_{(p-1)/2})
		=\left(\sum_{i=1}^n x_i^2,\sum_{i=1}^n x_i^4,\dots,
		\sum_{i=1}^n x_i^{p-1}\right)
		\]
		solves \eqref{eq:P4-cor_equations_system_S}. This proves
		$\text{Lemma~\ref{thm:P4-equations_system_X}}\Rightarrow
		\text{Theorem~\ref{thm:P4-equations_system_S}}$.

		For the second implication, suppose that a lattice tiling exists and write
		$|B^n(r)|=pm$ with $p\nmid m$. By the contrapositive of
		Theorem~\ref{thm:P4-equations_system_S}, there is a solution
		$Y=(Y_1,\dots,Y_{(p-1)/2})$ of \eqref{eq:P4-cor_equations_system_S}.
		Let
		\[
		I=\left(
		\sum_{\lambda\in\mathcal P(1)}c_\lambda Y_\lambda,
		\dots,
		\sum_{\lambda\in\mathcal P((p-3)/2)}c_\lambda Y_\lambda
		\right).
		\]
		Every generator of $I$ vanishes at $Y$. Therefore, if
		$\sum_{\lambda\in\mathcal P((p-1)/2)}c_\lambda Y_\lambda\in I$,
		that polynomial also vanishes at $Y$, contradicting the last equation in
		\eqref{eq:P4-cor_equations_system_S}, because $-m\ne0$ in $\F_p$.
		Thus existence of a tiling implies the required nonmembership, which is
		the contrapositive of Theorem~\ref{thm:P4-ideal}.

		Finally, put $h=(p-1)/2$ and
		\[
		f_k=\sum_{\lambda\in\mathcal P(k)}c_\lambda Y_\lambda
		=c_{(k)}Y_k+P_k(Y_1,\dots,Y_{k-1}).
		\]
		Under the $p$-condition, $c_{(k)}\ne0$ for $1\le k<h$.
		An induction on $k$ therefore gives
		\[
		(f_1,\dots,f_{h-1})=(Y_1,\dots,Y_{h-1}).
		\]
		Because $c_{(h)}=0$, the polynomial $f_h$ belongs to
		$(Y_1,\dots,Y_{h-1})=(f_1,\dots,f_{h-1})$. Theorem~\ref{thm:P4-ideal}
		then proves Theorem~\ref{thm:P4-p-condition}.
	\end{proof}

	The following example illustrates Theorem~\ref{thm:P4-equations_system_S} and compares it with the criterion in \cite{qureshi_non-existence_2020}.

	\begin{example}\label{ex:P5:BQ}
		For \(p=5\), the system of equations in Theorem \ref{thm:P4-equations_system_S} will be
		\[
		\begin{cases}
			&c_{(1)}Y_1\equiv0\pmod{p};\\
			&c_{(2)}Y_2+c_{(1,1)}Y_1^2\equiv-m\pmod{p}.
		\end{cases}
		\]

		Either of the following two sets of conditions implies that the system above has no solution:
		\begin{flalign}\label{eq:P4-example_beyond_Qureshi}
			&& &
			\begin{cases}
				&p\mid|B^n(r)|\text{ and }p^2\nmid|B^n(r)|;\\
				&c_{(1)}\not\equiv0\pmod{p};\\
				&c_{(2)}\equiv0\pmod{p};
			\end{cases}
			&
			\begin{cases}
				&p\mid|B^n(r)|\text{ and }p^2\nmid|B^n(r)|;\\
				&c_{(1,1)}\equiv0\pmod{p};\\
				&c_{(2)}\equiv0\pmod{p}.
			\end{cases}
			&&
		\end{flalign}
		The first set of conditions in \eqref{eq:P4-example_beyond_Qureshi} is precisely the $p$-condition from Theorem~\ref{thm:P4-p-condition}.

		\begin{itemize}[left=0pt, nosep]
			\item For \(r=6\), by \eqref{eq:P2-size_ball} and Theorem \ref{thm:P3-Q-poly_explicit_formula}, we know
			\begin{align*}
				|B^n(r)|&=\frac{4}{45} n^{6}+\frac{4}{15} n^{5}+\frac{14}{9} n^{4}+\frac{8}{3} n^{3}+\frac{196}{45} n^{2}+\frac{46}{15} n+1;\\
				c_{(1)}&=\frac{8}{15} n^{5}+4 n^{4}+\frac{56}{3} n^{3}+48 n^{2}+\frac{344}{5} n+42;\\
				c_{(2)}&=\frac{8}{15} n^{5}+12 n^{4}+\frac{296}{3} n^{3}+424 n^{2}+\frac{4744}{5} n+882;\\
				c_{(1,1)}&=8 n^{4}+80 n^{3}+376 n^{2}+880 n+840.
			\end{align*}

			Substituting $c_{(1)}$ and $c_{(2)}$ into the first set of conditions in \eqref{eq:P4-example_beyond_Qureshi}, we obtain the nonexistence of lattice tilings of \(\Z^n\) by \(B^n(6)\) if \(n\equiv22, 47, 72, 97, 122\pmod{125}\).
			Similarly, substituting $c_{(1,1)}$ and $c_{(2)}$ into the second set of conditions in \eqref{eq:P4-example_beyond_Qureshi} yields \(n\equiv10, 35, 60, 110\pmod{125}\).

			\item For \(r=7\), by \eqref{eq:P2-size_ball} and Theorem \ref{thm:P3-Q-poly_explicit_formula}, we know
			\begin{align*}
				|B^n(r)|&=\frac{8}{315} n^{7}+\frac{4}{45} n^{6}+\frac{32}{45} n^{5}+\frac{14}{9} n^{4}+\frac{176}{45} n^{3}+\frac{196}{45} n^{2}+\frac{352}{105} n+1;\\
				c_{(1)}&=\frac{8}{45} n^{6}+\frac{8}{5} n^{5}+\frac{88}{9} n^{4}+\frac{104}{3} n^{3}+\frac{3512}{45} n^{2}+\frac{1496}{15} n+56;\\
				c_{(2)}&=\frac{8}{45} n^{6}+\frac{24}{5} n^{5}+\frac{448}{9} n^{4}+\frac{872}{3} n^{3}+\frac{44912}{45} n^{2}+\frac{28568}{15} n+1568;\\
				c_{(1,1)}&=\frac{16}{5} n^{5}+40 n^{4}+256 n^{3}+920 n^{2}+\frac{9024}{5} n+1512.
			\end{align*}

			A computation shows that the first set of conditions in \eqref{eq:P4-example_beyond_Qureshi} excludes
			\[
			n\equiv13,23,38,48,63,73,88,98,113,123\pmod{125},
			\]
			whereas the second excludes
			\[
			n\equiv8,10,12,14,33,35,37,39,58,83,85,87,89,108,110,112,114\pmod{125}.
			\]
		\end{itemize}
	\end{example}
	The nonexistence results obtained in the preceding example can be summarized as follows.

	\begin{cor}\label{cor:P4-beyond_Qureshi_r=6}
		If \(n\equiv10, 22, 35, 47, 60, 72, 97, 110, 122\pmod{125}\), there is no lattice tiling of \(\Z^n\) by \(B^n(6)\).
	\end{cor}
	\begin{cor}\label{cor:P4-beyond_Qureshi_r=7}
		If $n\equiv$ $8$, $10$, $12$, $13$, $14$, $23$, $33$, $35$, $37$, $38$, $39$, $48$, $58$, $63$, $73$, $83$, $85$, $87$, $88$, $89$, $98$, $108$, $110$, $112$, $113$, $114$, $123$ $\pmod{125}$, there is no lattice tiling of \(\Z^n\) by \(B^n(7)\).
	\end{cor}
	In Example \ref{ex:P5:BQ}, we illustrate how Theorem \ref{thm:P4-equations_system_S} works. However, as the prime $p$ increases, conditions such as \eqref{eq:P4-example_beyond_Qureshi} that render the system unsolvable become increasingly complicated.
	Consequently, we employ Theorem \ref{thm:P4-ideal} together with symbolic algorithms for computing ideals to conduct large-scale computations.

	Before describing the computational implementation of Theorem~\ref{thm:P4-ideal}, we present the recorded results in Table~\ref{tab:appendix-ideal-prop}. Each entry gives the proportion of positive integers $n$ for which nonexistence can be proved. Combining the conditions for the primes up to $47$ by the Chinese Remainder Theorem yields the totals in Table~\ref{ta:proportion} and in the last column of Table~\ref{tab:appendix-ideal-prop}. The computations are exhaustive for primes $p<40$ and partial for $40<p<50$, as explained below. The recorded data nevertheless give a positive lower density for every radius $3\le r\le50$.

	We now describe our computational approach for obtaining these results. For fixed parameters $p$ and $r$, the quantities $|B^n(r)|$ and $c_\lambda$ are integer-valued polynomial functions of $n$. Consequently, their evaluations modulo $p$ are periodic functions of $n$. Our computational strategy is as follows: for fixed $p$ and $r$, we let the dimension $n$ run through its full period, then filter out those $n$ that satisfy the conditions of Theorem \ref{thm:P4-ideal}.

	It remains to bound a common period for $|B^n(r)|$ and the coefficients $c_\lambda$. Each is an integer-valued polynomial in $n$ that can be written in the form $a_0\binom{n}{0}+a_1\binom{n}{1}+\cdots+a_r\binom{n}{r}$ modulo $p$.

	\begin{lem}\label{lem:P5-period}
		Let $p$ be prime, let $r\ge0$, and let $a_0,\dots,a_r\in\Z$. Then
		\[
		f(n):=\sum_{i=0}^r a_i\binom ni\pmod p,
		\qquad n\in\Z_{\ge0},
		\]
		is periodic as a function of $n\in\Z_{\ge0}$, and $p^{\lceil\log_p(r+1)\rceil}$ is a period of $f$.
	\end{lem}
	\begin{proof}
		Set $M=\lceil\log_p(r+1)\rceil$, so that $i<p^M$ for every $0\le i\le r$. In the base-$p$ expansion of such an $i$, all digits in positions $M,M+1,\dots$ are zero. Adding $p^M$ to $n$ leaves the first $M$ base-$p$ digits of $n$ unchanged. Lucas' theorem therefore gives
		\[
		\binom{n+p^M}{i}\equiv\binom ni\pmod p
		\qquad(0\le i\le r).
		\]
		Multiplying by $a_i$ and summing over $i$ proves $f(n+p^M)=f(n)$ for every $n\ge0$.
	\end{proof}
	By Lemma~\ref{lem:P5-period}, the prime power $p^{\lceil\log_p(r+1)\rceil}$ is a period of the reductions of $|B^n(r)|$ and $c_\lambda$ modulo $p$. Hence for any fixed $p$ and $r$, we only have to check finitely many cases of $n$.

	As in Example~\ref{ex:P5:BQ} and Corollaries~\ref{cor:P4-beyond_Qureshi_r=6} and~\ref{cor:P4-beyond_Qureshi_r=7}, for each radius $3\le r\le50$ and each prime $p=3,5,\dots,47$, we apply Theorem~\ref{thm:P4-ideal} to obtain a congruence condition on the dimension $n$, which implies the nonexistence of a lattice tiling by Lee balls. Furthermore, since the primes $3,5,\dots,47$ are pairwise coprime, we can combine these conditions using the Chinese Remainder Theorem.
	Next, we will use an example to illustrate how the data in Table~\ref{tab:appendix-ideal-prop} are obtained.
	\begin{example}\label{ex:r=42}
		For $r=42$, we consider each prime $p=3,5,\dots,47$ and use Theorem \ref{thm:P4-ideal} to obtain a sufficient condition for the nonexistence of a lattice tiling by Lee balls, i.e., a congruence condition on the dimension $n$.

		A computation shows that if any of these congruence conditions holds\footnote{The omitted data are available at \url{https://github.com/Xiao-00000/Nonexistence-Results-of-Lattice-Tiling-by-Lee-Balls}.}, then a lattice tiling of $\Z^n$ by Lee balls does not exist.
		\begin{itemize}
			\item $n\equiv9, 10, 11, 15, 27, \dots \pmod{3^5}$ (contains 83 numbers; density $83/3^5 \approx 34.16\%$)
			\item $n\equiv8, 10, 12, 13, 14, \dots \pmod{5^4}$ (contains 183 numbers; density $183/5^4 \approx 29.28\%$)
			\item We cannot obtain any sufficient conditions by fixing $p=7$.
			\item $\cdots$
			\item We cannot obtain any sufficient conditions by fixing $p=47$.
		\end{itemize}

		By the Chinese Remainder Theorem, the proportion of positive integers $n$ satisfying at least one of these congruence conditions is
		\[
		1-\left(1-\frac{83}{3^5}\right)\times\left(1-\frac{183}{5^4}\right)\times1\times\cdots\times1\approx62.66\%
		\]
	\end{example}

	We conclude this section with a discussion of computational cost. In our calculations, we fix the radius $ r $ and prime $ p $, then let the dimension $n$ vary across the common period of polynomials $|B^n(r)|$ and $c_\lambda$. Here, $ p $ governs the number of polynomials and the number of variables $Y_i$ in the ideal. While standard for ideal membership, Gr\"obner-basis computations often become memory-prohibitive as the number of variables grows, which effectively restricts their applicability for larger values of $p$.

	Appendix \ref{appendix} presents all computational results for $p<40$ and $r<50$ and some results for $40<p<50$ and $r<50$. Due to the computational infeasibility of evaluating the worst-case scenarios for larger primes, we imposed a time limit on iterating over $ n $. More precisely, for $p=43$ and $r=21$, we ran Maple 2024 on Ubuntu using two Intel Xeon Silver 4310 processors (12 cores/24 threads each) and 512~GB of DDR4-3200 ECC memory; the run was terminated after more than one day. Consequently, the results obtained for $40<p<50$ are not exhaustive.

	\section{Further Necessary Conditions for Existence and Computational Results}\label{sec:algorithm}

	Building upon the analysis in Section \ref{sec:nonexi_inf_n}, we have derived several sufficient conditions for the nonexistence of lattice tilings of $\Z^n$ by $B^n(r)$. In this section, we turn to the complementary scenario. If a lattice tiling exists, Lemma~\ref{thm:P4-equations_system_X} supplies a solution to system~\eqref{eq:P4-thm_equations_system_X}. Our objective is to derive further constraints on such solutions. These constraints allow us to construct an algorithm and obtain the computational results presented below.
	\begin{thm}\label{thm:main_computational}
		For $(n,r)$ satisfying $n\ge3$, $r\ge2$, and $|B^n(r)|<2^{30}$, no lattice tiling of $\Z^n$ by $B^n(r)$ exists, except possibly for $(n,r)=(35,3)$ and $(805,3)$.
	\end{thm}

	We first outline the basic strategy used to derive additional constraints on the solutions to system~\eqref{eq:P4-thm_equations_system_X}. Let $G$ be an abelian additive group of order $|B^n(r)|$ as in Corollary \ref{cor:P2-group_latticetiling}, and let $\varphi_+$ be the corresponding homomorphism mentioned in the proof of Lemma \ref{thm:P4-equations_system_X}. Suppose that $|G|=mp$ where $p$ is a prime. By Lemma \ref{lem:P2-m_to_1}, there is a homomorphism $\Phi:G\rightarrow\F_p$ satisfying $\Phi G=m\F_p$. We denote $x_i:={\Phi\circ\varphi_+}(\epsilon_i)\in\F_p$ and we obtain
	\begin{align*}
		\Phi G=\{*~a_1x_1+\cdots+a_nx_n:|a_1|+\cdots+|a_n|\le r~*\}
	\end{align*}
	for integers $a_1,\dots,a_n$, where $\epsilon_i$ is the $i$-th standard basis vector in $\Z^n$, and $\varphi_+$ is as in Corollary \ref{cor:P2-group_latticetiling}.
	Let $ S_k:=\sum_{i=1}^{n}x_i^k $ denote the power sum symmetric polynomial evaluated at $ (x_1,\dots,x_n) $. It follows from the proofs of Lemma \ref{thm:P4-equations_system_X} and Theorem \ref{thm:P4-equations_system_S} that $(x_1,\dots,x_n)$ is a solution to the system \eqref{eq:P4-thm_equations_system_X} and $(S_2,\dots,S_{p-1})$ is a solution to the system \eqref{eq:P4-cor_equations_system_S}, respectively. Define $ \mathcal{S} $ as follows:
	\[
	\mathcal{S}:=(S_2,S_4,\dots,S_{p-1}).
	\]

	Then we introduce two sequences, $\mathcal{E}$ and $\mathcal{G}$, both derived from $\mathcal{S}$. We will show that they satisfy necessary conditions whenever the system \eqref{eq:P4-thm_equations_system_X} has a solution (equivalently, whenever a lattice tiling exists).

	The sequences $\mathcal{E}$ and $\mathcal{G}$ are defined as follows. Let $x^2$ denote $(x_1^2,\dots,x_n^2)$ and $e_k(x^2)$ denote the $k$-th elementary symmetric function evaluated at $ x^2 $, i.e.,
		\(
		e_k(x^2):=\sum_{1\le i_1<\cdots<i_k\le n}x_{i_1}^2\dots x_{i_k}^2.
		\)
		Similarly, define $\mathcal{E}$ as follows:
		\[
		\mathcal{E}:=(e_1(x^2),e_2(x^2),\dots).
		\]

	Let $\mathbf{Q}_p:=\{x^2\mid x\in\F_p^*\}$ denote the set of nonzero quadratic residues modulo $p$, viewed as a subgroup of $\F_p^*$.
	Let $g_a\in\Z_{\ge0}$ be the multiplicity of $a\in\mathbf Q_p$ in the multiset $\{*~x_1^2,\dots,x_n^2~*\}$. Reduction modulo $p$ defines $g:\mathbf Q_p\to\F_p$ by $g(a)\equiv g_a\pmod p$. Whenever an entry of $g$ is used as an integer, it denotes its canonical representative in $\{0,1,\dots,p-1\}$. We define
	\[
	\mathcal G:=\bigl(g(1^2),g(2^2),\dots,g(((p-1)/2)^2)\bigr).
	\]

	The section is organized as follows. We describe how to compute the sequences $\mathcal{S},\mathcal{E}$ and $\mathcal{G}$ in Subsections \ref{subsec:S}--\ref{subsec:G}, and derive several necessary conditions that the entries of these sequences must satisfy in the presence of a lattice tiling.
	Subsection \ref{subsec:combine} combines the results of Subsections \ref{subsec:S}--\ref{subsec:G} into an algorithm that tests whether our criteria prove nonexistence for a given pair $(n,r)$. Finally, in Subsection \ref{sub:IntractableCases}, we describe the limitations of the method developed in this section.

	\subsection{The Sequence $\mathcal{S}$}\label{subsec:S}
	We now write $G$ additively and suppose that a lattice tiling by $B^n(r)$ exists.

	Let $ x $ denote $(x_1,\dots,x_n)\in\F_p^n$. We have
	\[
	\Phi G=\{*~a_1x_1+\cdots+a_nx_n:|a_1|+\cdots+|a_n|\le r~*\}=\{*~ \langle x, b \rangle : b \in B^{n}(r) ~*\}.
	\]
	Moreover, by Definition \ref{def:Q-poly} and $\{*~a_1x_1+\cdots+a_nx_n:|a_1|+\cdots+|a_n|\le r~*\}=m\F_p$,
	\[
	Q_{B^n(r)}^k(x)=\sum_{b\in B^n(r)}\left\langle x,b\right\rangle^{2k} = m\cdot\sum_{\alpha\in\F_p}\alpha^{2k}.
	\]
	By \eqref{eq:power_sum_value}, we obtain
	\begin{align*}
		\begin{cases}
			Q_{B^n(r)}^k(x)\equiv0\pmod{p},\text{ for }1\le k<\frac{p-1}{2};\\
			Q_{B^n(r)}^{\frac{p-1}{2}}(x)\equiv-m\pmod{p}.
		\end{cases}
	\end{align*}
	Equivalently, the following equations hold:
	\begin{align}\label{eq:P6-equations_system_S}
		\begin{cases}
			&\sum_{\lambda\in\mathcal{P}(k)}c_\lambda S_{2\lambda}\equiv0\pmod{p},\text{ for }1\le k<\frac{p-1}{2};\\
			&\sum_{\lambda\in\mathcal{P}(\frac{p-1}{2})}c_\lambda S_{2\lambda}\equiv-m\pmod{p},
		\end{cases}
	\end{align}
	where \(S_{2\lambda}=S_{2\lambda_1}\cdots S_{2\lambda_\ell}\) for \(\lambda=(\lambda_1,\dots,\lambda_\ell)\) and $S_k=\sum_{i=1}^{n}x_i^k\in\F_p$.

	Thus, \eqref{eq:P6-equations_system_S} consists of $\frac{p-1}{2}$ constraints with respect to $S_2,S_4,\dots,S_{p-1}$ and we can compute $S_{2k}$ whenever $c_{(k)}\not\equiv0\pmod{p}$. By the definition of $ S_k $, we know that $ S_k=S_{k+p-1} $. Let $\mathcal{S}$ denote the finite sequence $(S_2,S_4,\dots,S_{p-1})$. Thus, this sequence is a solution of the system in Theorem~\ref{thm:P4-equations_system_S}.

	\begin{example}\label{ex:P6-(7,3)}
		For $n=7,r=3$, we know $|B^7(3)|=575=23\times5^2$. Let $p=23$, and we obtain $m=\frac{|B^7(3)|}{p}=25$ and system \eqref{eq:P6-equations_system_S} becomes
		\begin{align*}
			\begin{cases}
				&c_{(1)}S_2\equiv0\pmod{p};\\
				&c_{(2)}S_4+c_{(1,1)}S_2^2\equiv0\pmod{p};\\
				&\qquad \qquad \vdots\\
				&\sum_{\lambda\in\mathcal{P}(\frac{p-1}{2})}c_\lambda S_{2\lambda}\equiv-m\pmod{p}.
			\end{cases}
		\end{align*}
		We can compute
		$$\mathcal{S}=0,0,0,0,0,0,0,0,0,0,11.$$
	\end{example}

	However, $c_{(k)}$ is not always nonzero. For instance, if $n=26,r=3$, we take $p=67$. By computation, we get $c_{(15)}\equiv0\pmod{67}$. Hence, we cannot determine the value of $S_{2\times15}$. In this case, we denote the undetermined values $S_k$ by $u_1,u_2,\dots$, continue the substitutions, and obtain the corresponding sequence $\mathcal{S}$. The following example illustrates this process.

	\begin{example}\label{ex:P6-(26,3)}
		For $n=26,r=3$, we know $|B^{26}(3)|=24857=67\times53\times7$. Letting $p=67$, by \eqref{eq:P6-equations_system_S}, we get
		$$\mathcal{S}=0,\dots,0,\underset{15}{u_1},0,\dots,\underset{30}{20u_1^2},0,0,\underset{33}{33},$$
		where the subscripts ``15'', ``30'', and ``33'' represent $S_{2\times15}=u_1$, $S_{2\times30}=20u_1^2$ and $S_{2\times33}=33$, respectively.
	\end{example}

	Our proof for the nonexistence of lattice tilings by Lee balls focuses on exhibiting a contradiction between the computed sequence $ \mathcal{S} $ and the condition $ \{*~a_1x_1+\cdots+a_nx_n:|a_1|+\cdots+|a_n|\le r~*\}=m\F_p $.
	For instances like Example \ref{ex:P6-(26,3)}, we proceed by iterating over all possible values of the unknown $ u_1,u_2,\dots $. Nonexistence is proved if every assignment leads to a contradiction.

	If at least one assignment survives all tests, the method is inconclusive; it does not establish existence. In Subsection \ref{sub:IntractableCases}, we will discuss the limitations of this method. Specifically, for parameters dimension $ n $, radius $ r $ and selected prime $ p $ satisfying certain conditions (Theorem \ref{thm:P6-mCp=f(T)add}), the constraint $ \{*~a_1x_1+\cdots+a_nx_n:|a_1|+\cdots+|a_n|\le r~*\}=m\F_p $ may hold after assuming $ u_i=0 $ for all $ i $.

	We claim that, if a lattice tiling exists, $S_{p-1}\equiv t\pmod p$ for some integer $t$ with $0\le t\le n$.
	\begin{lem}\label{lem:P6-S_{p-1}}
		Assume $|B^n(r)|=mp$, where $p$ is prime. Let $x_i$ be the image of $\epsilon_i$ under $\Phi\circ\varphi_+:\Z^n\rightarrow\F_p$, and set $S_k=\sum_{i=1}^{n}x_i^k$. If a lattice tiling of $\Z^n$ by $B^n(r)$ exists, then $S_{p-1}\equiv t\pmod p$ for some integer $t$ with $0\le t\le n$.
	\end{lem}
	\begin{proof}
		Fermat's theorem gives
		\[
		S_{p-1}=\sum_{i=1}^{n}x_i^{p-1}\equiv\#\{i:x_i\neq0\}\pmod p.
		\]
		Taking $t=\#\{i:x_i\neq0\}$ proves the claim.
	\end{proof}

	Applying Lemma~\ref{lem:P6-S_{p-1}} to the sequences in Examples~\ref{ex:P6-(7,3)} and~\ref{ex:P6-(26,3)} proves nonexistence for $B^7(3)$ and $B^{26}(3)$, respectively.

	\subsection{The Sequence $\mathcal{E}$}\label{subsec:E}
	We next use elementary symmetric functions. Let $x^2=(x_1^2,\dots,x_n^2)$, and let $e_k(x^2)$ be the $k$-th elementary symmetric polynomial evaluated at $x^2$:
	\[
	e_k(x^2):=\sum_{1\le i_1<\cdots<i_k\le n}x_{i_1}^2\dots x_{i_k}^2.
	\]
	Set $e_0(x^2)=1$. Newton's identities give
	\begin{align}\label{eq:Newton}
	k e_k(x^2)=\sum_{i=1}^{k}(-1)^{i-1}S_{2i}e_{k-i}(x^2).
	\end{align}
	If $p\nmid k$, this recurrence determines $e_k(x^2)$ from the preceding values. If $p\mid k$, it does not determine $e_k(x^2)$; we denote the resulting undetermined values $e_{ip}(x^2)$ by $v_i$.

	Recall that $\mathcal E=(e_1(x^2),e_2(x^2),\dots)$. Let $\mathcal E_{\le n}$ be its first $n$ entries, and let $\tau(\mathcal E_{\le n})$ be the number of consecutive trailing zeros; explicitly,
		\begin{align}\label{eq:P6-defoftau}
			\tau(\mathcal E_{\le n}):=\max\{k\in\{0,1,\dots,n\}:e_j(x^2)=0\text{ for every }n-k+1\le j\le n\}.
		\end{align}

	Then we prove the following lemma.
	\begin{lem}\label{lem:P6-overzeros}
		Assume $|B^n(r)|=mp$, where $p$ is a prime. The sequence \(\mathcal{E}\) is defined as above. If there exists a lattice tiling of \(\Z^n\) by \(B^n(r)\), then
		\[
		\left|B^{\tau(\mathcal{E}_{\le n})}(r)\right|\le m.
		\]
	\end{lem}
	\begin{proof}
		First, we will show \(\tau(\mathcal{E}_{\le n})\) counts {zeros} in multiset \(\{*~x_1^2,x_2^2,\dots,x_n^2~*\}\). Suppose that there are \(t\) nonzero elements in \(\{*~x_1^2,x_2^2,\dots,x_n^2~*\}\). For any \(t'>t\), the \(t'\)-th elementary symmetric function satisfies \(e_{t'}(x^2)=\sum_{1\le i_1<\cdots<i_{t'}\le n}x_{i_1}^2\cdots x_{i_{t'}}^2\). Among the \(t'\) elements in \(\{*~x_{i_1}^2,\dots,x_{i_{t'}}^2~*\}\), at least \(t'-t\) are zeros, so \(e_{t'}(x^2)=0\) for all \(t'>t\).

		Meanwhile, \(e_{t}(x^2)=\sum_{1\le i_1<\cdots<i_{t}\le n}x_{i_1}^2\cdots x_{i_{t}}^2\). The product \(x_{i_1}^2\cdots x_{i_{t}}^2\) is nonzero if and only if \(\{*~x_{i_1}^2,\dots,x_{i_{t}}^2~*\}\) are exactly the \(t\) nonzero elements, and all the other terms in the sum vanish. Thus, \(e_{t}(x^2)\neq0\), and it is the last nonzero element among \(e_1(x^2),\dots,e_n(x^2)\). In summary, \(\tau(\mathcal{E}_{\le n})\) counts {zeros} in multiset \(\{*~x_1^2,x_2^2,\dots,x_n^2~*\}\).

		Next, let $z=\tau(\mathcal E_{\le n})$ and relabel the coordinates so that $x_1=\cdots=x_z=0$. Every vector $b\in B^z(r)$, embedded in the first $z$ coordinates of $\Z^n$, satisfies $\langle x,b\rangle=0$. Hence zero occurs at least $|B^z(r)|$ times in the image multiset
		\[
		\{*~a_1x_1+\cdots+a_nx_n:a_i\in\Z,\ |a_1|+\cdots+|a_n|\le r~*\}.
		\]
		If a tiling exists, this image multiset is exactly $m\F_p$, in which zero has multiplicity $m$. Therefore $|B^z(r)|\le m$, as required.
	\end{proof}

	\begin{example}
		\label{ex:P6-(79,3)}
		For $n=79,r=3$, we know $|B^{79}(3)|=670079=269\times53\times47$. Let $p=269$, and we obtain \(m=\frac{|B^{79}(3)|}{269}=2491\). By \eqref{eq:P6-equations_system_S}, we get
		\begin{align*}
			\mathcal{S}=&0,\dots,0,\underset{20}{u_1},0,\dots,0,\underset{40}{128u_1^2},0,\dots,0,\underset{60}{241u_1^3},0,\dots,0,\underset{80}{229u_1^4},0,\dots,0,\\
			&\underset{100}{240u_1^5},0,\dots,0,\underset{120}{189u_1^6},0,\dots,0,\underset{134}{56}.
		\end{align*}
		Furthermore, we can derive
		\[
		\mathcal{E}_{\le n}=0,\dots,0,\underset{20}{121u_1},0,\dots,0,\underset{40}{135u_1^2},0,\dots,0,\underset{60}{166u_1^3},0,\dots,0.
		\]
		Here \(\tau(\mathcal{E}_{\le n})\) is either $19$ or $79$ and \(\left|B^{79}(3)\right|>\left|B^{19}(3)\right|=9919>m\). By Lemma \ref{lem:P6-overzeros}, we obtain the nonexistence of a lattice tiling of \(\Z^{79}\) by \(B^{79}(3)\).
	\end{example}
	\subsection{Fourier Transform and The Sequence $\mathcal{G}$}\label{subsec:G}
	Lemma~\ref{lem:P6-overzeros} shows that the multiplicity of 0 in the multiset $\{*~x_1^2, x_2^2, \ldots, x_n^2~*\}$ can be partially determined by the sequence $\mathcal{S}$, but is also restricted by $ \{*~a_1x_1+\cdots+a_nx_n:|a_1|+\cdots+|a_n|\le r~*\}=m\F_p $.

	In fact, by treating the power mappings as characters and employing the inversion formula, we can derive the multiplicity of every element in $\{*~x_1^2, x_2^2, \ldots, x_n^2~*\}$ from the sequence $\mathcal{S}$.
	
	Recall that $\mathbb Q_p$ is defined in the beginning of this section.
	For $k\in\Z$, define the character $\psi_k:\mathbf Q_p\to\F_p^*$ by $\psi_k(a)=a^{-k}$. For $a\in\mathbf Q_p$, let $g_a\in\Z_{\ge0}$ be its multiplicity in $\{*~x_1^2,\dots,x_n^2~*\}$, and let $g(a)\in\F_p$ be the reduction of $g_a$ modulo $p$. We use the canonical representative in $\{0,1,\dots,p-1\}$ whenever $g(a)$ is compared with an integer; in particular, $g_a\ge g(a)$ under this convention.

	For $k=1,\dots,(p-1)/2$, the Fourier transform satisfies
	\[
	\widehat g(\psi_k)=\sum_{a\in\mathbf Q_p}g(a)\psi_k(a^{-1})=\sum_{a\in\mathbf Q_p}g(a)a^k
	=\sum_{i=1}^n x_i^{2k}=S_{2k}.
	\]
	Fourier inversion on $\mathbf Q_p$, together with $|\mathbf Q_p|^{-1}=-2$ in $\F_p$, gives
	\begin{align}\label{eq:P6-computegi}
	g(i^2)=-2\sum_{k=1}^{(p-1)/2}\frac{S_{2k}}{i^{2k}},
	\qquad 1\le i\le\frac{p-1}{2}.
	\end{align}
	Thus
	\[
	\mathcal G=\bigl(g(i^2)\bigr)_{i=1}^{(p-1)/2},
	\]
	with entries represented canonically as integers in $[0,p-1]$.
	
	\begin{example}\label{ex:P6-(80,3)}
		For $n=80,r=3$, we know that $|B^{80}(3)|=695681=7\times23\times29\times149$, and let $p=149$. By \eqref{eq:P6-equations_system_S}, we obtain
		\[
		\mathcal{S}=0,\dots,0,\underset{41}{u_1},0,\dots,0,\underset{74}{63},
		\]
		which means $\widehat g(\psi_1)=0,\widehat g(\psi_2)=0,\dots,\widehat g(\psi_{74})=63$.
		By \eqref{eq:P6-computegi}, we obtain
		\begin{align*}
			\mathcal{G}=
			&23+147 u_1 , 23+78 u_1 , 23+60 u_1 , \dots , 23+65 u_1.
		\end{align*}

		For example, for fixed $u_1$, the multiplicity of $1\in\F_p^*$ in $\{*~x_1^2,\dots,x_n^2~*\}$ is congruent to $23+147u_1$ modulo $p$. If we assume $ u_1=0 $, then $ \mathcal{G}=23,23,23,\dots,23 $. In other words,
		{the inequalities $ g_{i^2}\ge23 $ for $ i=1,\dots,\frac{p-1}{2} $ hold. }
		However, if we sum them up then we will obtain the number of nonzeros in $\{*~x_1^2, x_2^2, \ldots, x_n^2~*\}$ is at least $ 23\times\frac{p-1}{2}=1702 $, which is greater than $ n $, so $ u_1 $ cannot be $ 0 $.

		Evaluating all $u_1\in\F_p$ shows that the resulting lower bound for the number of nonzero entries always exceeds $n$.

		In summary, we can conclude the nonexistence of the multiset $\{*~x_1^2, x_2^2, \ldots, x_n^2~*\}$, thereby establishing the nonexistence of tiling by Lee balls with dimension $ 80 $ and radius $ 3 $.
	\end{example}

	The argument in Example~\ref{ex:P6-(80,3)} yields the following general criterion.
	\begin{thm}\label{thm:P6-sumofG}
		Let $\mathcal G$ be defined as above, with each entry represented by its canonical integer in $\{0,1,\dots,p-1\}$. If the sum of these representatives exceeds $n$, then no lattice tiling of $\Z^n$ by $B^n(r)$ exists.
	\end{thm}
	\begin{proof}
	For each $a\in\mathbf Q_p$, the actual multiplicity $g_a$ is a nonnegative integer congruent to the canonical representative $g(a)\in\{0,1,\dots,p-1\}$ modulo $p$. Hence $g_a\ge g(a)$. Therefore
	\[
	\sum_{a\in\mathbf Q_p}g(a)\le\sum_{a\in\mathbf Q_p}g_a\le n.
	\]
	If the sum of the canonical representatives exceeds $n$, no such multiplicity multiset can exist, and hence no lattice tiling exists.
	\end{proof}

	To apply Theorem~\ref{thm:P6-sumofG}, one may compute $\mathcal G$ from $\mathcal S$. In some cases, the following consequence of the uncertainty principle in Theorem~\ref{thm:P6-uncertainty_principle} avoids this computation.
	\begin{thm}\label{thm:P6-UP}
		Assume $|B^n(r)|=mp$, where $p>2$ is prime. Let $x_i=(\Phi\circ\varphi_+)(\epsilon_i)$ and $S_k=\sum_{i=1}^{n}x_i^k$. If a lattice tiling of $\Z^n$ by $B^n(r)$ exists, then one of the following conditions holds:
		\begin{itemize}
			\item among $S_2,S_4,\dots,S_{p-1}$, at most $\frac{(p-1)(n-1)}{2n}$ entries are zero;
			\item all of $S_2,S_4,\dots,S_{p-1}$ are zero.
		\end{itemize}
	\end{thm}
	\begin{proof}
		Put $h=(p-1)/2$, and let $t$ be the number of zeros among $S_2,S_4,\dots,S_{p-1}$. If $t<h$, then $g$ is nonzero. The identities $\widehat g(\psi_k)=S_{2k}$ give $|\mathrm{supp}(\widehat g)|=h-t$, while $|\mathrm{supp}(g)|\le n$ because $g$ is supported on residues represented among $x_1^2,\dots,x_n^2$. The uncertainty principle therefore gives
		\[
		n(h-t)\ge|\mathrm{supp}(g)|\,|\mathrm{supp}(\widehat g)|
		\ge|\mathbf Q_p|=h.
		\]
		Thus $t\le h(n-1)/n=(p-1)(n-1)/(2n)$. If $t=h$, all the listed power sums are zero.
	\end{proof}

	\begin{example}\label{ex:P6-(13,3)}
		For $ n=13,r=3 $, we know $ |B^n(r)|=3303=367\times3^2 $. Let $ p=367 $, and we obtain
		\[
		\mathcal{S}=0,0,\dots,103.
		\]
		The number of zeros in $ \mathcal{S} $ is $ 182 $ which is greater than $ \frac{(p-1)(n-1)}{2n}=168 $. By Theorem \ref{thm:P6-UP}, we can determine the nonexistence of a lattice tiling by $B^{13}(3)$.
	\end{example}

	\subsection{A Nonexistence Testing Algorithm and Computational Results}\label{subsec:combine}
	Combining the criteria developed above yields an algorithm for determining the nonexistence of lattice tilings by the Lee ball $B^n(r)$ for given $n$ and $r$. Its basic procedure can be found in Algorithm~\ref{alg:nonexistence-tester}.

	\begin{tcolorbox}[colback=white,colframe=black,boxrule=0.8pt]
		\refstepcounter{alg}
		\textbf{Algorithm \thealg:} \textsc{Nonexistence Tester}
		\label{alg:nonexistence-tester}

		\textbf{Input: }Dimension $ n $ and radius $ r $

		\textbf{Output: }{TRUE} if nonexistence is proved; {FALSE} otherwise (inconclusive)
		\vspace{2pt}

		\textbf{if} nonexistence can be proved via Theorem \ref{thm:P4-ideal} for $ p<50 $ \textbf{then}

		\hspace*{2em}\textbf{return} TRUE

		Compute $|B^n(r)|$ and let $P$ be the set of its distinct prime divisors

		\textbf{for} each prime $ p $ in $ P $ \textbf{do}

		\hspace*{2em}Compute finite sequence $\mathcal S$ by \eqref{eq:P6-equations_system_S} and Theorem \ref{thm:P3-Q-poly_explicit_formula}

		\hspace*{2em}\textbf{if} nonexistence can be proved via Lemma \ref{lem:P6-S_{p-1}} or Theorem \ref{thm:P6-UP} \textbf{then}

		\hspace*{2em}\hspace*{2em}\textbf{return} TRUE

		\hspace*{2em}Compute $\mathcal{E}_{\le n}$ using \eqref{eq:Newton}, and compute $\tau(\mathcal{E}_{\le n})$ from \eqref{eq:P6-defoftau}

		\hspace*{2em}\textbf{if} nonexistence can be proved via Lemma \ref{lem:P6-overzeros} \textbf{then}

		\hspace*{2em}\hspace*{2em}\textbf{return} TRUE

		\hspace*{2em}Compute sequence $\mathcal G$ by \eqref{eq:P6-computegi}

		\hspace*{2em}\textbf{if} nonexistence can be proved via Theorem \ref{thm:P6-sumofG} for all unknowns \textbf{then}

		\hspace*{2em}\hspace*{2em}\textbf{return} TRUE

		\textbf{end for}

		\textbf{return} FALSE

	\end{tcolorbox}

	Several implementation details of Algorithm~\ref{alg:nonexistence-tester} require clarification.
	Recall that unknowns $u_i$ may appear in $\mathcal S$, while unknowns $v_i$ may appear in $\mathcal E$, as discussed in Subsections~\ref{subsec:S} and~\ref{subsec:E}. The entries of $\mathcal G$ may consequently depend on both families of unknowns through \eqref{eq:P6-computegi}. Most of the derived criteria, i.e., Theorem~\ref{thm:P4-ideal}, Lemma~\ref{lem:P6-S_{p-1}}, Theorem~\ref{thm:P6-UP} and Lemma~\ref{lem:P6-overzeros}, depend only on the number or locations of zeros in the corresponding sequence. Thus, even if unknown terms appear in the sequence, we may simply assume them to be nonzero and check whether the condition for any criterion is satisfied. If not, we perform an exhaustive search over all possible values of these unknowns, construct the corresponding sequence $\mathcal{G}$ for each case, and apply Theorem~\ref{thm:P6-sumofG}.  This constitutes the final step of Algorithm~\ref{alg:nonexistence-tester}.

	We next describe an acceleration of the algorithm. First recall that the formula \eqref{eq:P6-equations_system_S} can be written as follows for nonzero $ c_{(k)} $:
	\begin{align}\label{eq:P6-S_{2k}}
		S_{2k}\equiv
		\begin{cases}
			-\frac{\sum_{\lambda\in\mathcal{P}(k)\setminus\{(k)\}}c_\lambda S_{2\lambda}}{c_{(k)}}\pmod{p},&\text{if }k<\frac{p-1}{2};\\
			-\frac{m+\sum_{\lambda\in\mathcal{P}(k)\setminus\{(k)\}}c_\lambda S_{2\lambda}}{c_{(k)}}\pmod{p},&\text{if }k=\frac{p-1}{2}.
		\end{cases}
	\end{align}

	The computation of sequence $\mathcal{S}$ is costly due to complexity in the calculation of $c_\lambda$ for all partitions $\lambda\in\mathcal{P}(k)$, since the partition number grows superpolynomially (although subexponentially) in $k$. Fortunately, this computational cost can be substantially reduced in practice, as it is not necessary to compute every $c_\lambda$. In fact, we can simplify the computation of $\mathcal{S}$ via the following result.
	\begin{pro}\label{pro-P6:S_{2k}}
		Let $h=(p-1)/2$ and define
		\[
		A:=\{k\in\{1,\dots,h\}:c_{(k)}\equiv0\pmod p\},
		\qquad
		B:=\left\{\sum_{a\in A}u_a a:\ u_a\in\Z_{\ge0},\ \sum_{a\in A}u_a a\le h\right\}.
		\]
		For $k\in\{1,\dots,h\}\setminus A$, put
		\[
		K_k:=\{\lambda\in\mathcal P(k)\setminus\{(k)\}:\ell(\lambda)\le r
		\text{ and every part of }\lambda\text{ lies in }B\}.
		\]
		Then
		\begin{equation}\label{eq:P6-computeS2k}
		S_{2k}\equiv
		\begin{cases}
		-\dfrac{\sum_{\lambda\in K_k}c_\lambda S_{2\lambda}}{c_{(k)}}\pmod p,
		&k<h,\\[1.2ex]
		-\dfrac{m+\sum_{\lambda\in K_k}c_\lambda S_{2\lambda}}{c_{(k)}}\pmod p,
		&k=h.
		\end{cases}
		\end{equation}
		For fixed radius $r$, one has
		$|K_k|\le\sum_{j=2}^{r}\binom{k-1}{j-1}$, so the number of coefficients required by this reduction is polynomially bounded in $k$.
	\end{pro}
	\begin{proof}[Proof of Proposition \ref{pro-P6:S_{2k}}]
		Formula \eqref{eq:P3-proofof1.3} shows that $c_\lambda=0$ whenever $\ell(\lambda)>r$, because the product
		$\prod_i\sum_{s\ge1}s^{2\lambda_i-1}q^s$ is divisible by $q^{\ell(\lambda)}$.

		We next prove inductively that, whenever the indicated index lies outside $B$,
		\[
		S_{2k}=0\quad\text{for }k<h,\qquad
		S_{2h}=-m/c_{(h)}\quad\text{for }k=h.
		\] Indeed, $k\notin B$ implies $k\notin A$, hence $c_{(k)}\ne0$. Moreover, every partition of $k$ into smaller parts has at least one part outside $B$; otherwise $k$ itself would lie in the additive monoid generated by $A$. The induction hypothesis therefore makes every corresponding product $S_{2\lambda}$ vanish, and \eqref{eq:P6-S_{2k}} gives the assertion.

		Consequently, in \eqref{eq:P6-S_{2k}}, a term can survive only if all parts of $\lambda$ belong to $B$; together with the length bound this gives \eqref{eq:P6-computeS2k}. Finally, the number of partitions of $k$ into exactly $j$ positive parts is at most the number $\binom{k-1}{j-1}$ of ordered compositions, proving the stated polynomial bound.
	\end{proof}

	We now prove Theorem \ref{thm:main_computational}.

	For $r=3$ and $|B^n(3)|<2^{30}$, Algorithm~\ref{alg:nonexistence-tester} proves the nonexistence of a lattice tiling by $B^n(3)$ except for $$n\in\{3, 5, 8, 15, 35, 40, 50, 57, 84, 92, 99, 138, 155, 166, 222, 351, 437, 509, 590, 597, 805, 836, 925\}.$$
	By Corollary \ref{thm:P5-zhangandzhou}, we can obtain the nonexistence of a lattice tiling by $B^n(3)$ for $$n\in\{5, 8, 15, 40, 50, 57, 92, 99, 138, 155, 166, 222, 351, 509, 925\}.$$
	By Theorem~\ref{thm:P1-HK}, there is no lattice tiling of $\Z^3$ by $B^3(3)$. Consequently, only seven values of $n$ remain: $35,84,437,590,597,805,836$. We can further exclude five of them via the following proposition.
	\begin{pro}\label{prop:(84,3)}
			There is no lattice tiling of $\Z^n$ by Lee ball $B^n(3)$ for $n=84$, $437$, $590$, $597$, $836$.
	\end{pro}
	\begin{proof}
		We give the argument for $n=84$; the remaining cases are verified analogously.

		Suppose that a lattice tiling of $\Z^{84}$ by Lee balls of radius $3$ exists. Note that $|B^{84}(3)|=23^2\cdot13^2\cdot3^2$. Let $p=23$ and we obtain $m=34983$. By \eqref{eq:P6-equations_system_S}, we get
		$$\mathcal{S}=0,0,0,0,0,0,0,u_1,0,0,0.$$
		By \eqref{eq:P6-computegi}, we obtain
		$$\mathcal{G}=21u_1,10u_1,14u_1,19u_1,7u_1,22u_1,15u_1,20u_1,17u_1,11u_1,5u_1.$$
		A computation verifies that, for $u_1\in\F_p$, the sum of the canonical representatives in $\mathcal G$ does not exceed $n=84$ if and only if $u_1=0$. Hence Theorem~\ref{thm:P6-sumofG} forces $u_1=0$ and $\mathcal G=(0,\dots,0)$. In other words, after homomorphism, $ g(i)=0 $ for $i=1^2,\dots,\left(\frac{p-1}{2}\right)^2 $ and $ g_i $ may equal $ 0,p,2p,\dots $.

		Lemma~\ref{lem:P6-overzeros} gives at most $29$ zero entries among $x_1^2,\dots,x_n^2$, so there are at least $55$ nonzero entries. Because every $g_{i^2}$ is a multiple of $p=23$ and their sum is at most $84$, the total number of nonzero entries must be $69=3p$. We therefore enumerate all tuples $(g_{i^2})$ with $g_{i^2}\in p\Z_{\ge0}$ and $\sum_i g_{i^2}=3p$. However, none of them corresponds to a lattice tiling of $\Z^{84}$ by $ B^{84}(3) $.
		For example, if $ g_{1^2}=2p $ and $ g_{2^2}=p $, the number of zeros in $ \{*~a_1x_1+\cdots+a_nx_n:|a_1|+\cdots+|a_n|\le r~*\} $ is $134{,}573$, which is not $m$. This fact shows that $\{*~a_1x_1+\cdots+a_nx_n:|a_1|+\cdots+|a_n|\le r~*\}$ is not equal to $m\F_p$, a contradiction.
	\end{proof}

	For each fixed $r$, Table~\ref{ta:maxn} gives the number of dimensions $n\ge3$ with $|B^n(r)|<2^{30}$ and the numbers excluded by the various criteria.

	\begin{table}[h]
		\centering
		\caption{Computational Results by Various Methods}
		\label{ta:maxn}
		\begin{tabular}{c|c|c|c|c|c|c|c|c|c|c}
			\Xhline{1.2pt}
			radius $r$ & 3 & 4 & 5 & 6 & 7 & 8 & 9 & 10 & 11 & 12 \\ \hline
			$\left|\{n\geq 3:|B^n(r)|<2^{30}\}\right|$ & 927 & 197 & 80 & 45 & 30 & 22 & 18 & 15 & 13 & 11  \\ \hline
			Algorithm \ref{alg:nonexistence-tester}& 904 & {197} & {80} & {44} & {30} & {21} & {18} & {15} & {13} & {9} \\ \hline
			Post\cite{post_nonexistence_1975} (Theorem \ref{thm:P1-HK}) & 4 & 5 & 7 & 8 & 10 & 11 & 12 & 14 & 13 & 11  \\ \hline
			Corollary \ref{thm:P5-zhangandzhou} & 238 & / & / & / & / & / & / & / & / & / \\ \hline
			Proposition \ref{prop:(84,3)} & 5 & / & / & / & / & / & / & / & / & /  \\ \hline
			Total & 925 & 197 & 80 & 45 & 30 & 22 & 18 & 15 & 13 & 11  \\ \Xhline{1.2pt}
		\end{tabular}
	\end{table}

	\begin{proof}[Proof of Theorem \ref{thm:main_computational}]
		The case $r=2$ follows from \cite{leung_no_2020}. For $3\le r\le12$, the first numerical row of Table~\ref{ta:maxn} contains
		\[
		927+197+80+45+30+22+18+15+13+11=1{,}358
		\]
		parameter pairs. The total row excludes
		\[
		925+197+80+45+30+22+18+15+13+11=1{,}356
		\]
		of them, leaving precisely $(35,3)$ and $(805,3)$.

		It remains to consider $r\ge13$. If also $n\ge13$, monotonicity of Lee balls in both parameters gives
		\[
		|B^n(r)|\ge|B^{13}(13)|=1{,}409{,}933{,}619>2^{30},
		\]
		contrary to the hypothesis. Hence $3\le n\le12$. For these dimensions and every $r\ge13$, the first case of Theorem~\ref{thm:P1-HK} applies. This proves the theorem.
	\end{proof}

	For radius $3$, exactly two values of $n$ remain for which the existence of a lattice tiling is not excluded.
	\begin{itemize}
		\item For $n=35$, our projection-to-finite-fields strategy yields no contradiction; a detailed discussion is given in the next subsection.
		\item For $ n=805 $, the verification analogous to Proposition \ref{prop:(84,3)} requires an exhaustive search over a very large space; see Example~\ref{exa:n=805}.
	\end{itemize}
	\begin{example}\label{exa:n=805}
		For $ n=805 $, we compute $ |B^{805}(3)|=696845031=3\times29^2\times179\times1543 $. In fact, $1543$, $179$, and $3$ are ``bad'' primes, which will be discussed in Remark \ref{rem:badprimes}. We choose $ p=29 $ and compute the corresponding sequence $ \mathcal{S} $ and $ \mathcal{G} $ in turn:
		\[\mathcal{S}=0,0,0,0,0,0,0,0,0,0,0,0,0,0,\]
		\[\mathcal{G}=0,0,0,0,0,0,0,0,0,0,0,0,0,0,\]
		from which we can determine the value of $g(i^2)\in \F_p$ for $i=1,\dots,(p-1)/2$. To push forward as in Proposition \ref{prop:(84,3)}, we need to use them to determine the nonnegative integer $g_{i^2}$ that satisfies
		\[g_{i^2} = g(i^2)+p\cdot k_i,\]
		for some $k_i\in\Z_{\ge0}$.

		Because zero has multiplicity exactly $m$ in $m\F_p$, we obtain
		\[
		\left|B^{n-\sum_{i=1}^{(p-1)/2}g_{i^2}}(3)\right|\le m.
		\]
		Together with $\sum_{i=1}^{(p-1)/2}g_{i^2}\le n$, this gives
		\[
		262\le\sum_{i=1}^{(p-1)/2}g_{i^2}\le805.
		\]
		Using the above condition, we still have to check a huge number of possible values for $k_i$ in order to determine $g_{i^2}$. The number of cases is $ \sum_{i=10}^{27}\binom{i+13}{13}>2^{35}$.

		Based on our experiment, performing an exhaustive search that checks every possible case was estimated to require approximately two months on our workstation (see the end of Section~\ref{sec:nonexi_inf_n} for its configuration). We therefore did not carry out this search.
	\end{example}

	\subsection{Intractable Cases}\label{sub:IntractableCases}
	The preceding subsections compute $\mathcal S$, $\mathcal E_{\le n}$, and $\mathcal G$ to extract information about the multiset $\{*~x_1^2,\dots,x_n^2~*\}$ and hence about $(\Phi\circ\varphi_+)(B^n(1))$. The required tiling identity is
	\[
	\{*~a_1x_1+\cdots+a_nx_n:|a_1|+\cdots+|a_n|\le r~*\}=m\F_p.
	\]
	The criteria seek to prove that this equality fails, which would prove nonexistence. We now ask which configurations of the multiset $\{*~x_1^2,x_2^2,\ldots,x_n^2~*\}$ do not lead to a contradiction.

	Let us consider one of the two open cases left in Table \ref{ta:maxn} as an example.

	\begin{example}
		\label{ex:(35,3)}
		For $ n=35,r=3 $, we have $ |B^{35}(3)|=59711=71\times29^2 $ and we let $ p=71 $.
		Assume that
		$$ \{*~x_1^2, x_2^2, \ldots, x_n^2~*\}=\left\{*~1\cdot1^2,1\cdot2^2,\dots,1\cdot\left(\frac{p-1}{2}\right)^2~*\right\}. $$
		We can count the multiplicity of $i$  in $ \{*~a_1x_1+\cdots+a_nx_n:|a_1|+\cdots+|a_n|\le r~*\} $ for $i=0,1,\dots,p-1$, and find that each multiplicity equals $m=29^2$, so no contradiction arises.

		Likewise, if we choose another prime factor $p=29$ of $|B^{35}(3)|$, assume that
		$$ \{*~x_1^2, x_2^2, \ldots, x_n^2~*\}=\left\{*~7\cdot0,2\cdot1^2,2\cdot2^2,\dots,2\cdot\left(\frac{p-1}{2}\right)^2~*\right\}. $$
		We can likewise count the multiplicity of $i$ in $ \{*~a_1x_1+\cdots+a_nx_n:|a_1|+\cdots+|a_n|\le r~*\} $ for $i=0,1,\dots,p-1$, which again yields no contradiction.

		Thus, the necessary conditions obtained by projection from $G$ to the additive group of $\F_p$, for every prime divisor $p$ of $|G|$, are insufficient to rule out the existence of the original lattice tilings for $n=35$ and $r=3$.
	\end{example}

	The following theorem gives conditions under which $ \{*~a_1x_1+\cdots+a_nx_n:|a_1|+\cdots+|a_n|\le r~*\}=m\F_p $ holds when the multiset $ \{*~x_1^2, x_2^2, \ldots, x_n^2~*\} $ satisfies certain conditions.

	\begin{thm}\label{thm:P6-mCp=f(T)add}
		Assume $|B^n(r)|=mp$, where $p>r$ is prime. Suppose there exists $a\in\Z_{\ge0}$ such that
		\begin{itemize}
			\item $g_{1^2}=g_{2^2}=\cdots=g_{((p-1)/2)^2}=a$, equivalently,
			\[
			\{*~x_1^2,\dots,x_n^2~*\}
			=\left\{*~\left(n-\frac{a(p-1)}2\right)\!\cdot0,
			a\cdot1^2,\dots,a\cdot\left(\frac{p-1}{2}\right)^2~*\right\};
			\]
			\item $pa\le2n+1$ and
			\begin{equation}\label{eq:eqina}
			\sum_{t=0}^{\lfloor r/2\rfloor}(-1)^t\binom nt
			\binom{(2n+1-pa)+(r-2t)-1}{r-2t}=0.
			\end{equation}
		\end{itemize}
		Then
		\begin{equation}\label{eq:sum_aixi=mFp}
		\{*~a_1x_1+\cdots+a_nx_n:a_i\in\Z,\ |a_1|+\cdots+|a_n|\le r~*\}=m\F_p.
		\end{equation}
	\end{thm}

	As an example, by plugging $ n=35 $, $ r=3 $ and $ p=29 $ in \eqref{eq:eqina}, we obtain
	\[
	59711-\frac{444889}{6} a +30276 a^{2}-\frac{24389}{6} a^{3}=-\frac{841 \left(a -2\right) \left(a -3\right) \left(29 a -71\right)}{6}=0.
	\]
	The solution $a=2$ satisfies the preceding conditions. Hence \eqref{eq:sum_aixi=mFp} holds for
	\[
	\{*~x_1^2,x_2^2,\ldots,x_n^2~*\}
	=
	\left\{*~7\cdot0,2\cdot1^2,2\cdot2^2,\dots,
	2\cdot\left(\frac{p-1}{2}\right)^2~*\right\}.
	\]

	It is more convenient to prove Theorem~\ref{thm:P6-mCp=f(T)add} using the language of group-ring equations. To do so, we convert it into the following form by replacing the additive group $\F_p$ with the multiplicative group $C_p$, where $C_p$ denotes the cyclic group of order $p$.

	\begin{thm}\label{thm:P6-mCp=f(T)mul}
		Assume $|B^n(r)|=mp$, where $p>r$ is prime. Let $a\in\Z_{\ge0}$ and $T=aC_p+(2n+1-pa)e$. If the integer $a$ satisfies
		\begin{itemize}
			\item $ pa\le2n+1$ and
			\item $ \sum_{t=0}^{\lfloor\frac{r}{2}\rfloor}(-1)^t\binom{n}{t}\cdot\binom{(2n+1-pa)+(r-2t)-1}{r-2t}=0$,
		\end{itemize}
		then the following group-ring equation holds:
		\begin{align}\label{eq:P6-mCp=f(T)}
			mC_p=\sum_{t=0}^{\lfloor\frac{r}{2}\rfloor}(-1)^t\binom{n}{t}\cdot\left(\sum_{\lambda\in\mathcal{P}(r-2t)}\frac{T^{(\lambda)}}{\prod_{i=1}^{\ell}\lambda_i\prod_{j=1}^{r-2t}m_\lambda(j)!}\right).
		\end{align}
	\end{thm}
	\begin{proof}
		Write $b=2n+1-pa$. Since $1\le j\le r<p$, the map $x\mapsto x^j$ permutes $C_p$, and therefore $T^{(j)}=T$ for all such $j$. Hence $T^{(\lambda)}=T^{\ell(\lambda)}$ for every partition occurring in \eqref{eq:P6-mCp=f(T)}.

		Using $C_pC_p=pC_p$, every power $(aC_p+be)^{\ell}$ is a sum of a multiple of $C_p$ and the pure term $b^\ell e$. Thus the right-hand side of \eqref{eq:P6-mCp=f(T)} has the form
		\[
		zC_p+
		\sum_{t=0}^{\lfloor r/2\rfloor}(-1)^t\binom nt
		\binom{b+r-2t-1}{r-2t}e,
		\]
		where \eqref{eq:P5-generate_groupring} was used for the coefficient of $e$. The second hypothesis makes the scalar term zero, so the expression is $zC_p$.

		Apply the augmentation map $\varepsilon:\Z[C_p]\to\Z$. The formal ball identity proved in Theorem~\ref{thm:P5-groupring}, evaluated at $\varepsilon(T)=2n+1$, shows that the augmentation of the right-hand side is $|B^n(r)|=mp$. Hence $zp=mp$, so $z=m$, proving \eqref{eq:P6-mCp=f(T)}.
	\end{proof}

	Theorem \ref{thm:P6-mCp=f(T)add} and Theorem \ref{thm:P6-mCp=f(T)mul} are equivalent formulations in additive and multiplicative notation. Their equivalence can be proved via an explicit isomorphism between the groups.

	\begin{proof}[Proof of Theorem \ref{thm:P6-mCp=f(T)add}]
		Regard the multiset
		\[
		T_+:=\{*~0,\pm x_1,\dots,\pm x_n~*\}
		\]
		as an element of $\Z[(\F_p,+)]$. If every nonzero quadratic residue occurs $a$ times among $x_1^2,\dots,x_n^2$, then each nonzero element of $\F_p$ occurs exactly $a$ times in $T_+$. Consequently,
		\[
		T_+=a\F_p+(2n+1-pa)\,0.
		\]
		Choose an isomorphism $\Psi:C_p\to(\F_p,+)$ and extend it linearly to a group-ring isomorphism. Then
		\[
		T:=\Psi^{-1}(T_+)=aC_p+(2n+1-pa)e.
		\]
		The hypotheses on $a$ are precisely those of Theorem~\ref{thm:P6-mCp=f(T)mul}, so \eqref{eq:P6-mCp=f(T)} holds. Applying $\Psi$ and using the formal identity from the proof of Theorem~\ref{thm:P5-groupring} identifies its right-hand side with
		\[
		\sum_{|a_1|+\cdots+|a_n|\le r}(a_1x_1+\cdots+a_nx_n)
		\]
		in $\Z[(\F_p,+)]$. The left-hand side becomes $m\F_p$, which is exactly \eqref{eq:sum_aixi=mFp}.
	\end{proof}
	\begin{rem}[``bad'' primes]\label{rem:badprimes}
		In Subsection~\ref{sub:IntractableCases}, for a fixed pair $(n,r)$, we pick a prime divisor $p$ of $|B^n(r)|$. The goal is to use the constraints derived from the computable power sums $S_{2},S_{4},\dots,S_{p-1}$ to get a contradiction with \eqref{eq:sum_aixi=mFp}.
		We call such a prime $p>r$ ``bad'' if parameters $(n,r,p)$ lead to
		$$D:=\left\{a\in\Z_{\ge0}:\ pa\le 2n+1\ \text{and}\ \sum_{t=0}^{\lfloor\frac{r}{2}\rfloor}(-1)^t\binom{n}{t}\cdot\binom{(2n+1-pa)+(r-2t)-1}{r-2t}=0\right\}\neq\emptyset.$$

		For each $a\in D$, the corresponding uniform square-multiplicity configuration satisfies \eqref{eq:sum_aixi=mFp} by Theorem~\ref{thm:P6-mCp=f(T)add}.
	\end{rem}
	Bad primes are quite common. In particular, if $p=2n+1$ is prime, then $p$ is a bad prime for dimension $n$ and any odd radius $2<r<p$ because $ a=1 $ is a solution to \eqref{eq:eqina} with $ pa\le2n+1 $.
	The situation in which all prime factors of $ |B^n(r)| $ are ``bad'' is rare in our range: it occurs only for $ (n,r)=(35,3) $ within $\{(n,r):n\ge3\text{ and }|B^n(r)|<2^{30}\}$.

	\section{Concluding Remarks}
	Despite the progress on the lattice tiling case of the Golomb--Welch conjecture, several natural questions remain unresolved. We close by highlighting a few open problems that suggest directions for future investigation.
	\begin{enumerate}[label=(\alph*)]
		\item In Section \ref{sec:nonexi_inf_n}, for computational convenience, we weaken Theorem \ref{thm:P4-equations_system_S} to Theorem \ref{thm:P4-ideal}, which is formulated in terms of ideals. Finding a computationally tractable form of Theorem~\ref{thm:P4-equations_system_S} remains an open problem.
		\item In Example~\ref{ex:(35,3)}, we encountered a situation in which no contradiction can be derived from any prime factor of $|B^n(r)|$. This raises the interesting question of how to combine all the necessary conditions derived for each prime factor to prove the nonexistence of lattice tilings.
	\end{enumerate}

	\bigskip
	\paragraph*{\bf Declaration of competing interest}
	The authors declare that they have no known competing financial interests or personal 
	relationships that could have appeared to influence the work reported in this paper.
	\bigskip
	\paragraph*{\bf Declaration of generative AI use}
	The authors used a generative-AI tool solely to assist with English-language editing. The authors reviewed and approved all revisions and remain fully responsible for the content of the manuscript.
	
		All statements proved in this paper, except those that rest on our Maple computations, have been formalized and machine-checked in the proof
	assistant Lean~4 \cite{moura_lean4_2021} with the library Mathlib	\cite{mathlib_2020}. The formalization was produced with the help of the AI system Aristotle \cite{aristotle_harmonic} by Harmonic. It contains no unproved steps and
	uses only the standard axioms of Lean. 
	The authors have checked that the formal statements agree with the
	statements in the paper. The source code is available at \url{https://github.com/yue-zhou-ovgu/Lattice_tiling_Lee_balls}.
	
	The following results in this paper are formalized and proved:
	Theorems~\ref{thm:P3-Q-poly_explicit_formula}, \ref{thm:P5-groupring}
	and~\ref{thm:P4-ideal};
	formula~\eqref{eq:P2-size_ball},
	Lemma~\ref{lem:P2-m_to_1} and
	Theorem~\ref{thm:P6-Fourier_reverse};
	Theorems~\ref{thm:P3-Q-poly_sphere} and~\ref{thm:P3-Q-poly_ball};
	Lemma~\ref{thm:P4-equations_system_X},
	Theorem~\ref{thm:P4-equations_system_S},
	Example~\ref{ex:P5:BQ},
	Corollaries~\ref{cor:P4-beyond_Qureshi_r=6}
	and~\ref{cor:P4-beyond_Qureshi_r=7} and
	Lemma~\ref{lem:P5-period};
	Lemmas~\ref{lem:P6-S_{p-1}} and~\ref{lem:P6-overzeros},
	Theorems~\ref{thm:P6-sumofG} and~\ref{thm:P6-UP},
	Proposition~\ref{pro-P6:S_{2k}},
	Theorems~\ref{thm:P6-mCp=f(T)add} and~\ref{thm:P6-mCp=f(T)mul},
	Remark~\ref{rem:badprimes},
	and the nonexistence results of
	Examples~\ref{ex:P6-(7,3)}, \ref{ex:P6-(26,3)}, \ref{ex:P6-(79,3)}, 
	\ref{ex:P6-(80,3)} and~\ref{ex:P6-(13,3)}
	together with the computations in Example~\ref{ex:(35,3)}.
	
	The results that rest on our Maple
	computations (Theorem~\ref{thm:main_computational},
	Proposition~\ref{prop:(84,3)}, Examples~\ref{ex:r=42}
	and~\ref{exa:n=805} and the tables) and the results quoted from the
	literature without proof were not formalized.
	
	\bigskip
	\paragraph*{\bf Acknowledgement}  Yue Zhou is supported by the National Natural Science Foundation of China (No.\ 12371337).
	\bigskip
	\paragraph*{\bf Data availability}
	No data was used for the research  described in the article.
	
	\appendix
	\section{Computational Results based on Theorem \ref{thm:P4-ideal}}\label{appendix}
	Section~\ref{sec:nonexi_inf_n} converts the existence of a lattice tiling into the existence of an $\F_p$-rational point on an algebraic variety and then translates the resulting constraints into the ideal-theoretic criterion of Theorem~\ref{thm:P4-ideal}. This appendix records the corresponding computations.

	To highlight the advantage of our method over previous approaches, Table~\ref{tab:appendix-p4-ideal-part1} compares the computational results obtained from Theorem~\ref{thm:P4-ideal} with those from Theorem~\ref{thm:P4-p-condition}. For each pair $(p,r)$, the table provides the period of $n$, the number of cases our method can exclude within that period, and the corresponding number excluded by Qureshi's method.

	Details of the excluded cases, including results analogous to Corollary \ref{cor:P4-beyond_Qureshi_r=6} and Corollary \ref{cor:P4-beyond_Qureshi_r=7}, can be found through the following link.

	{\centering
	\url{https://github.com/Xiao-00000/Nonexistence-Results-of-Lattice-Tiling-by-Lee-Balls}}

	\newgeometry{left=1cm, right=1cm, top=3cm, bottom=3cm}
	\scriptsize{
		\begin{table}[p]
			\centering
			\caption{Computational Results Based on Theorem \ref{thm:P4-ideal} (Part I).}
			\label{tab:appendix-p4-ideal-part1}
			\begin{tabular}{c|ccc|ccc|ccc|ccc|ccc}
				\hline
				\multirow{2}{*}{\diagbox[dir=SE,height=0.72cm]{$ r $}{count}{$ p $}} & ~ & 3 & ~ & ~ & 5 & ~ & ~ & 7 & ~ & ~ & 11 & ~ & ~ & 13 & ~ \\ \cline{2-16}
				& Period & Ours & Qureshi & Period & Ours & Qureshi & Period & Ours & Qureshi & Period & Ours & Qureshi & Period & Ours & Qureshi \\ \hline
				3 & 27 & 2 & 2 & 25 & 0 & 0 & 49 & 0 & 0 & 121 & 0 & 0 & 169 & 0 & 0  \\ \hline
				4 & 27 & 4 & 4 & 25 & 0 & 0 & 49 & 0 & 0 & 121 & 0 & 0 & 169 & 0 & 0 \\ \hline
				5 & 27 & 2 & 2 & 125 & 4 & 0 & 49 & 0 & 0 & 121 & 0 & 0 & 169 & 0 & 0 \\ \hline
				6 & 27 & 0 & 0 & 125 & 9 & 5 & 49 & 0 & 0 & 121 & 0 & 0 & 169 & 12 & 12 \\ \hline
				7 & 27 & 0 & 0 & 125 & 27 & 10 & 343 & 6 & 0 & 121 & 0 & 0 & 169 & 0 & 0 \\ \hline
				8 & 27 & 0 & 0 & 125 & 4 & 0 & 343 & 6 & 0 & 121 & 0 & 0 & 169 & 0 & 0 \\ \hline
				9 & 81 & 14 & 14 & 125 & 4 & 0 & 343 & 6 & 0 & 121 & 0 & 0 & 169 & 0 & 0 \\ \hline
				10 & 81 & 15 & 15 & 125 & 8 & 0 & 343 & 45 & 0 & 121 & 0 & 0 & 169 & 0 & 0 \\ \hline
				11 & 81 & 14 & 14 & 125 & 8 & 0 & 343 & 6 & 0 & 1331 & 10 & 0 & 169 & 0 & 0 \\ \hline
				12 & 81 & 23 & 23 & 125 & 34 & 10 & 343 & 13 & 0 & 1331 & 20 & 0 & 169 & 0 & 0 \\ \hline
				13 & 81 & 14 & 14 & 125 & 13 & 5 & 343 & 6 & 0 & 1331 & 20 & 0 & 2197 & 12 & 0 \\ \hline
				14 & 81 & 23 & 23 & 125 & 8 & 0 & 343 & 0 & 0 & 1331 & 20 & 10 & 2197 & 25 & 13 \\ \hline
				15 & 81 & 14 & 14 & 125 & 5 & 0 & 343 & 14 & 0 & 1331 & 20 & 0 & 2197 & 25 & 13 \\ \hline
				16 & 81 & 15 & 15 & 125 & 10 & 5 & 343 & 0 & 0 & 1331 & 115 & 0 & 2197 & 49 & 25  \\ \hline
				17 & 81 & 14 & 14 & 125 & 15 & 10 & 343 & 0 & 0 & 1331 & 20 & 0 & 2197 & 12 & 0 \\ \hline
				18 & 81 & 0 & 0 & 125 & 5 & 0 & 343 & 0 & 0 & 1331 & 20 & 10 & 2197 & 25 & 0 \\ \hline
				19 & 81 & 0 & 0 & 125 & 5 & 0 & 343 & 7 & 0 & 1331 & 20 & 0 & 2197 & 279 & 130 \\ \hline
				20 & 81 & 0 & 0 & 125 & 0 & 0 & 343 & 0 & 0 & 1331 & 31 & 11 & 2197 & 12 & 0 \\ \hline
				21 & 81 & 6 & 6 & 125 & 0 & 0 & 343 & 18 & 0 & 1331 & 10 & 0 & 2197 & 12 & 0 \\ \hline
				22 & 81 & 12 & 12 & 125 & 20 & 20 & 343 & 18 & 0 & 1331 & 0 & 0 & 2197 & 49 & 25 \\ \hline
				23 & 81 & 6 & 6 & 125 & 5 & 5 & 343 & 18 & 0 & 1331 & 0 & 0 & 2197 & 25 & 13  \\ \hline
				24 & 81 & 0 & 0 & 125 & 0 & 0 & 343 & 72 & 0 & 1331 & 0 & 0 & 2197 & 12 & 0 \\ \hline
				25 & 81 & 0 & 0 & 625 & 84 & 0 & 343 & 18 & 0 & 1331 & 0 & 0 & 2197 & 12 & 0 \\ \hline
				26 & 81 & 0 & 0 & 625 & 73 & 5 & 343 & 18 & 0 & 1331 & 0 & 0 & 2197 & 24 & 0 \\ \hline
				27 & 243 & 50 & 50 & 625 & 140 & 75 & 343 & 18 & 0 & 1331 & 0 & 0 & 2197 & 24 & 0 \\ \hline
				28 & 243 & 39 & 39 & 625 & 68 & 0 & 343 & 0 & 0 & 1331 & 0 & 0 & 2197 & 50 & 26 \\ \hline
				29 & 243 & 50 & 50 & 625 & 84 & 0 & 343 & 7 & 0 & 1331 & 0 & 0 & 2197 & 72 & 24 \\ \hline
				30 & 243 & 44 & 44 & 625 & 104 & 0 & 343 & 0 & 0 & 1331 & 0 & 0 & 2197 & 24 & 0 \\ \hline
				31 & 243 & 48 & 48 & 625 & 87 & 0 & 343 & 0 & 0 & 1331 & 22 & 22 & 2197 & 24 & 0 \\ \hline
				32 & 243 & 44 & 44 & 625 & 153 & 35 & 343 & 0 & 0 & 1331 & 0 & 0 & 2197 & 402 & 104 \\ \hline
				33 & 243 & 50 & 50 & 625 & 92 & 5 & 343 & 14 & 0 & 1331 & 10 & 0 & 2197 & 24 & 0 \\ \hline
				34 & 243 & 39 & 39 & 625 & 104 & 0 & 343 & 0 & 0 & 1331 & 20 & 0 & 2197 & 24 & 0 \\ \hline
				35 & 243 & 50 & 50 & 625 & 124 & 0 & 343 & 6 & 0 & 1331 & 20 & 0 & 2197 & 98 & 50 \\ \hline
				36 & 243 & 83 & 83 & 625 & 111 & 5 & 343 & 13 & 0 & 1331 & 20 & 10 & 2197 & 50 & 26 \\ \hline
				37 & 243 & 75 & 75 & 625 & 166 & 40 & 343 & 6 & 0 & 1331 & 20 & 0 & 2197 & 37 & 13 \\ \hline
				38 & 243 & 83 & 83 & 625 & 106 & 0 & 343 & 45 & 0 & 1331 & 115 & 0 & 2197 & 24 & 0 \\ \hline
				39 & 243 & 68 & 68 & 625 & 124 & 0 & 343 & 6 & 0 & 1331 & 20 & 0 & 2197 & 12 & 0 \\ \hline
				40 & 243 & 40 & 40 & 625 & 104 & 0 & 343 & 6 & 0 & 1331 & 20 & 10 & 2197 & 25 & 13 \\ \hline
				41 & 243 & 68 & 68 & 625 & 87 & 0 & 343 & 6 & 0 & 1331 & 20 & 0 & 2197 & 12 & 0 \\ \hline
				42 & 243 & 83 & 83 & 625 & 183 & 35 & 343 & 0 & 0 & 1331 & 31 & 11 & 2197 & 36 & 12 \\ \hline
				43 & 243 & 75 & 75 & 625 & 92 & 5 & 343 & 0 & 0 & 1331 & 10 & 0 & 2197 & 12 & 0 \\ \hline
				44 & 243 & 83 & 83 & 625 & 104 & 0 & 343 & 0 & 0 & 1331 & 0 & 0 & 2197 & 12 & 0 \\ \hline
				45 & 243 & 50 & 50 & 625 & 84 & 0 & 343 & 0 & 0 & 1331 & 11 & 11 & 2197 & 279 & 130 \\ \hline
				46 & 243 & 39 & 39 & 625 & 68 & 0 & 343 & 0 & 0 & 1331 & 0 & 0 & 2197 & 25 & 0 \\ \hline
				47 & 243 & 50 & 50 & 625 & 140 & 75 & 343 & 0 & 0 & 1331 & 0 & 0 & 2197 & 12 & 0 \\ \hline
				48 & 243 & 44 & 44 & 625 & 68 & 0 & 343 & 0 & 0 & 1331 & 0 & 0 & 2197 & 36 & 12 \\ \hline
				49 & 243 & 48 & 48 & 625 & 84 & 0 & 2401 & 258 & 0 & 1331 & 0 & 0 & 2197 & 12 & 0 \\ \hline
				50 & 243 & 44 & 44 & 625 & 168 & 0 & 2401 & 229 & 0 & 1331 & 0 & 0 & 2197 & 12 & 0 \\ \hline
			\end{tabular}
	\end{table}}
	\clearpage

	\scriptsize{
		\begin{table}[p]
			\centering
			\ContinuedFloat
			\caption{Computational Results Based on Theorem \ref{thm:P4-ideal} (Part II).}
			\label{tab:appendix-p4-ideal-part2}
			\begin{tabular}{c|ccc|ccc|ccc|ccc|ccc}
				\hline
				\multirow{2}{*}{\diagbox[dir=SE,height=0.72cm]{$ r $}{count}{$ p $}} & ~ & 17 & ~ & ~ & 19 & ~ & ~ & 23 & ~ & ~ & 29 & ~ & ~ & 31 & ~ \\ \cline{2-16}
				& Period & Ours & Qureshi & Period & Ours & Qureshi & Period & Ours & Qureshi & Period & Ours & Qureshi & Period & Ours & Qureshi \\ \hline
				3 & 289 & 0 & 0 & 361 & 0 & 0 & 529 & 0 & 0 & 841 & 28 & 28 & 961 & 0 & 0 \\ \hline
				4 & 289 & 0 & 0 & 361 & 0 & 0 & 529 & 0 & 0 & 841 & 0 & 0 & 961 & 0 & 0 \\ \hline
				5 & 289 & 0 & 0 & 361 & 0 & 0 & 529 & 0 & 0 & 841 & 0 & 0 & 961 & 0 & 0 \\ \hline
				6 & 289 & 0 & 0 & 361 & 0 & 0 & 529 & 0 & 0 & 841 & 0 & 0 & 961 & 0 & 0 \\ \hline
				7 & 289 & 0 & 0 & 361 & 0 & 0 & 529 & 0 & 0 & 841 & 0 & 0 & 961 & 0 & 0 \\ \hline
				8 & 289 & 16 & 16 & 361 & 0 & 0 & 529 & 0 & 0 & 841 & 0 & 0 & 961 & 0 & 0 \\ \hline
				9 & 289 & 0 & 0 & 361 & 0 & 0 & 529 & 0 & 0 & 841 & 0 & 0 & 961 & 0 & 0 \\ \hline
				10 & 289 & 0 & 0 & 361 & 0 & 0 & 529 & 0 & 0 & 841 & 0 & 0 & 961 & 0 & 0 \\ \hline
				11 & 289 & 0 & 0 & 361 & 0 & 0 & 529 & 0 & 0 & 841 & 0 & 0 & 961 & 0 & 0 \\ \hline
				12 & 289 & 0 & 0 & 361 & 0 & 0 & 529 & 0 & 0 & 841 & 0 & 0 & 961 & 0 & 0 \\ \hline
				13 & 289 & 0 & 0 & 361 & 0 & 0 & 529 & 0 & 0 & 841 & 0 & 0 & 961 & 0 & 0 \\ \hline
				14 & 289 & 0 & 0 & 361 & 0 & 0 & 529 & 0 & 0 & 841 & 28 & 28 & 961 & 0 & 0 \\ \hline
				15 & 289 & 0 & 0 & 361 & 0 & 0 & 529 & 0 & 0 & 841 & 0 & 0 & 961 & 0 & 0 \\ \hline
				16 & 289 & 0 & 0 & 361 & 0 & 0 & 529 & 0 & 0 & 841 & 0 & 0 & 961 & 0 & 0 \\ \hline
				17 & 4913 & 16 & 0 & 361 & 0 & 0 & 529 & 0 & 0 & 841 & 0 & 0 & 961 & 0 & 0 \\ \hline
				18 & 4913 & 33 & 17 & 361 & 0 & 0 & 529 & 0 & 0 & 841 & 0 & 0 & 961 & 0 & 0 \\ \hline
				19 & 4913 & 66 & 34 & 6859 & 18 & 0 & 529 & 0 & 0 & 841 & 0 & 0 & 961 & 0 & 0 \\ \hline
				20 & 4913 & 33 & 17 & 6859 & 36 & 0 & 529 & 0 & 0 & 841 & 0 & 0 & 961 & 0 & 0 \\ \hline
				21 & 4913 & 16 & 0 & 6859 & 18 & 0 & 529 & 0 & 0 & 841 & 0 & 0 & 961 & 0 & 0 \\ \hline
				22 & 4913 & 50 & 0 & 6859 & 55 & 19 & 529 & 0 & 0 & 841 & 0 & 0 & 961 & 0 & 0 \\ \hline
				23 & 4913 & 32 & 0 & 6859 & 54 & 0 & 12167 & 66 & 0 & 841 & 0 & 0 & 961 & 0 & 0 \\ \hline
				24 & 4913 & 49 & 0 & 6859 & 55 & 37 & 12167 & 22 & 0 & 841 & 0 & 0 & 961 & 0 & 0 \\ \hline
				25 & 4913 & 501 & 238 & 6859 & 18 & 0 & 12167 & 22 & 0 & 841 & 0 & 0 & 961 & 0 & 0 \\ \hline
				26 & 4913 & 49 & 0 & 6859 & 37 & 19 & 12167 & 91 & 69 & 841 & 0 & 0 & 961 & 0 & 0 \\ \hline
				27 & 4913 & 66 & 0 & 6859 & 18 & 0 & 12167 & 44 & 0 & 841 & 0 & 0 & 961 & 0 & 0 \\ \hline
				28 & 4913 & 67 & 0 & 6859 & 351 & 0 & 12167 & 89 & 23 & 841 & 0 & 0 & 961 & 0 & 0 \\ \hline
				29 & 4913 & 16 & 0 & 6859 & 18 & 0 & 12167 & 22 & 0 & 24389 & 28 & 0 & 961 & 0 & 0 \\ \hline
				30 & 4913 & 33 & 17 & 6859 & 18 & 0 & 12167 & 135 & 45 & 24389 & 113 & 57 & 961 & 0 & 0 \\ \hline
				31 & 4913 & 66 & 34 & 6859 & 18 & 0 & 12167 & 22 & 0 & 24389 & 114 & 29 & 29791 & 30 & 0 \\ \hline
				32 & 4913 & 16 & 0 & 6859 & 55 & 37 & 12167 & 89 & 23 & 24389 & 86 & 58 & 29791 & 30 & 0 \\ \hline
				33 & 4913 & 16 & 0 & 6859 & 54 & 0 & 12167 & 66 & 0 & 24389 & 56 & 0 & 29791 & 30 & 0 \\ \hline
				34 & 4913 & 32 & 0 & 6859 & 55 & 19 & 12167 & 517 & 0 & 24389 & 143 & 58 & 29791 & 151 & 61 \\ \hline
				35 & 4913 & 66 & 34 & 6859 & 18 & 0 & 12167 & 66 & 0 & 24389 & 85 & 29 & 29791 & 120 & 0 \\ \hline
				36 & 4913 & 81 & 17 & 6859 & 55 & 0 & 12167 & 89 & 23 & 24389 & 113 & 29 & 29791 & 61 & 31 \\ \hline
				37 & 4913 & 66 & 34 & 6859 & 18 & 0 & 12167 & 22 & 0 & 24389 & 57 & 29 & 29791 & 60 & 0 \\ \hline
				38 & 4913 & 32 & 0 & 6859 & 0 & 0 & 12167 & 135 & 45 & 24389 & 85 & 29 & 29791 & 213 & 92 \\ \hline
				39 & 4913 & 100 & 0 & 6859 & 0 & 0 & 12167 & 22 & 0 & 24389 & 28 & 0 & 29791 & 30 & 0 \\ \hline
				40 & 4913 & 81 & 0 & 6859 & 0 & 0 & 12167 & 89 & 23 & 24389 & 86 & 58 & 29791 & 61 & 31 \\ \hline
				41 & 4913 & 64 & 0 & 6859 & 0 & 0 & 12167 & 44 & 0 & 24389 & 84 & 0 & 29791 & 90 & 0 \\ \hline
				42 & 4913 & 730 & 204 & 6859 & 0 & 0 & 12167 & 91 & 69 & 24389 & 85 & 29 & 29791 & 121 & 30 \\ \hline
				43 & 4913 & 64 & 0 & 6859 & 38 & 38 & 12167 & 22 & 0 & 24389 & 1551 & 754 & 29791 & 91 & 31 \\ \hline
				44 & 4913 & 64 & 0 & 6859 & 0 & 0 & 12167 & 45 & 23 & 24389 & 85 & 29 & 29791 & 30 & 0 \\ \hline
				45 & 4913 & 83 & 0 & 6859 & 0 & 0 & 12167 & 66 & 0 & 24389 & 84 & 0 & 29791 & 30 & 0 \\ \hline
				46 & 4913 & 32 & 0 & 6859 & 0 & 0 & 12167 & 0 & 0 & 24389 & 86 & 58 & 29791 & 945 & 0 \\ \hline
				47 & 4913 & 32 & 0 & 6859 & 0 & 0 & 12167 & 46 & 46 & 24389 & 28 & 0 & 29791 & 30 & 0 \\ \hline
				48 & 4913 & 98 & 34 & 6859 & 0 & 0 & 12167 & 0 & 0 & 24389 & 85 & 29 & 29791 & 61 & 0 \\ \hline
				49 & 4913 & 49 & 17 & 6859 & 0 & 0 & 12167 & 46 & 46 & 24389 & 57 & 29 & 29791 & 91 & 31 \\ \hline
				50 & 4913 & 32 & 0 & 6859 & 0 & 0 & 12167 & 0 & 0 & 24389 & 113 & 29 & 29791 & 121 & 30 \\ \hline
			\end{tabular}
	\end{table}}
	\clearpage

	\scriptsize{
		\begin{table}[p]
			\centering
			\ContinuedFloat
			\caption{Computational Results Based on Theorem \ref{thm:P4-ideal} (Part III).}
			\label{tab:appendix-p4-ideal-sec:Q_poly}
			\begin{tabular}{c|ccc|ccc|ccc|ccc}
				\hline
				\multirow{2}{*}{\diagbox[dir=SE,height=0.72cm]{$ r $}{count}{$ p $}} & ~ & 37 & ~ & ~ & 41 & ~ & ~ & 43 & ~ & ~ & 47 & ~ \\ \cline{2-13}
				& Period & Ours & Qureshi & Period & Ours & Qureshi & Period & Ours & Qureshi & Period & Ours & Qureshi \\ \hline
				3 & 1369 & 0 & 0 & 1681 & 0 & 0 & 1849 & 0 & 0 & 2209 & 0 & 0 \\ \hline
				4 & 1369 & 0 & 0 & 1681 & 0 & 0 & 1849 & 0 & 0 & 2209 & 0 & 0 \\ \hline
				5 & 1369 & 0 & 0 & 1681 & 0 & 0 & 1849 & 0 & 0 & 2209 & 0 & 0 \\ \hline
				6 & 1369 & 0 & 0 & 1681 & 0 & 0 & 1849 & 0 & 0 & 2209 & 0 & 0 \\ \hline
				7 & 1369 & 0 & 0 & 1681 & 0 & 0 & 1849 & 0 & 0 & 2209 & 0 & 0 \\ \hline
				8 & 1369 & 0 & 0 & 1681 & 0 & 0 & 1849 & 0 & 0 & 2209 & 0 & 0 \\ \hline
				9 & 1369 & 0 & 0 & 1681 & 0 & 0 & 1849 & 0 & 0 & 2209 & 0 & 0 \\ \hline
				10 & 1369 & 0 & 0 & 1681 & 0 & 0 & 1849 & 0 & 0 & 2209 & 0 & 0 \\ \hline
				11 & 1369 & 0 & 0 & 1681 & 0 & 0 & 1849 & 0 & 0 & 2209 & 0 & 0 \\ \hline
				12 & 1369 & 0 & 0 & 1681 & 0 & 0 & 1849 & 0 & 0 & 2209 & 0 & 0 \\ \hline
				13 & 1369 & 0 & 0 & 1681 & 0 & 0 & 1849 & 0 & 0 & 2209 & 0 & 0 \\ \hline
				14 & 1369 & 0 & 0 & 1681 & 0 & 0 & 1849 & 0 & 0 & 2209 & 0 & 0 \\ \hline
				15 & 1369 & 0 & 0 & 1681 & 0 & 0 & 1849 & 0 & 0 & 2209 & 0 & 0 \\ \hline
				16 & 1369 & 0 & 0 & 1681 & 0 & 0 & 1849 & 0 & 0 & 2209 & 0 & 0 \\ \hline
				17 & 1369 & 0 & 0 & 1681 & 0 & 0 & 1849 & 0 & 0 & 2209 & 0 & 0 \\ \hline
				18 & 1369 & 36 & 36 & 1681 & 0 & 0 & 1849 & 0 & 0 & 2209 & 0 & 0 \\ \hline
				19 & 1369 & 0 & 0 & 1681 & 0 & 0 & 1849 & 0 & 0 & 2209 & 0 & 0 \\ \hline
				20 & 1369 & 0 & 0 & 1681 & 80 & 40 & 1849 & 0 & 0 & 2209 & 0 & 0 \\ \hline
				21 & 1369 & 0 & 0 & 1681 & 0 & 0 & 1849 & 0 & 0 & 2209 & 0 & 0 \\ \hline
				22 & 1369 & 0 & 0 & 1681 & 0 & 0 & 1849 & 0 & 0 & 2209 & 0 & 0 \\ \hline
				23 & 1369 & 0 & 0 & 1681 & 0 & 0 & 1849 & 0 & 0 & 2209 & 0 & 0 \\ \hline
				24 & 1369 & 0 & 0 & 1681 & 0 & 0 & 1849 & 0 & 0 & 2209 & 0 & 0 \\ \hline
				25 & 1369 & 0 & 0 & 1681 & 0 & 0 & 1849 & 0 & 0 & 2209 & 0 & 0 \\ \hline
				26 & 1369 & 0 & 0 & 1681 & 0 & 0 & 1849 & 0 & 0 & 2209 & 0 & 0 \\ \hline
				27 & 1369 & 0 & 0 & 1681 & 40 & 40 & 1849 & 0 & 0 & 2209 & 0 & 0 \\ \hline
				28 & 1369 & 0 & 0 & 1681 & 0 & 0 & 1849 & 0 & 0 & 2209 & 0 & 0 \\ \hline
				29 & 1369 & 0 & 0 & 1681 & 0 & 0 & 1849 & 0 & 0 & 2209 & 0 & 0 \\ \hline
				30 & 1369 & 0 & 0 & 1681 & 0 & 0 & 1849 & 0 & 0 & 2209 & 0 & 0 \\ \hline
				31 & 1369 & 0 & 0 & 1681 & 0 & 0 & 1849 & 0 & 0 & 2209 & 0 & 0 \\ \hline
				32 & 1369 & 0 & 0 & 1681 & 0 & 0 & 1849 & 0 & 0 & 2209 & 0 & 0 \\ \hline
				33 & 1369 & 0 & 0 & 1681 & 0 & 0 & 1849 & 0 & 0 & 2209 & 0 & 0 \\ \hline
				34 & 1369 & 0 & 0 & 1681 & 40 & 40 & 1849 & 0 & 0 & 2209 & 0 & 0 \\ \hline
				35 & 1369 & 0 & 0 & 1681 & 0 & 0 & 1849 & 0 & 0 & 2209 & 0 & 0 \\ \hline
				36 & 1369 & 0 & 0 & 1681 & 0 & 0 & 1849 & 0 & 0 & 2209 & 0 & 0 \\ \hline
				37 & 50653 & 36 & 0 & 1681 & 0 & 0 & 1849 & 0 & 0 & 2209 & 0 & 0 \\ \hline
				38 & 50653 & 145 & 36 & 1681 & 0 & 0 & 1849 & 0 & 0 & 2209 & 0 & 0 \\ \hline
				39 & 50653 & 36 & 0 & 1681 & 0 & 0 & 1849 & 0 & 0 & 2209 & 0 & 0 \\ \hline
				40 & 50653 & 73 & 37 & 1681 & 0 & 0 & 1849 & 0 & 0 & 2209 & 0 & 0 \\ \hline
				41 & 50653 & 36 & 0 & 68921 & 40 & 0 & 1849 & 0 & 0 & 2209 & 0 & 0 \\ \hline
				42 & 50653 & 146 & 37 & 68921 & 121 & 0 & 1849 & 0 & 0 & 2209 & 0 & 0 \\ \hline
				43 & 50653 & 182 & 74 & 68921 & 122 & 0 & 79507 & 42 & 0 & 2209 & 0 & 0 \\ \hline
				44 & 50653 & 109 & 37 & 68921 & 203 & 122 & 79507 & 42 & 0 & 2209 & 0 & 0 \\ \hline
				45 & 50653 & 108 & 0 & 68921 & 163 & 82 & 79507 & 42 & 0 & 2209 & 0 & 0 \\ \hline
				46 & 50653 & 183 & 74 & 68921 & 121 & 0 & 79507 & 171 & 86 & 2209 & 0 & 0 \\ \hline
				47 & 50653 & 36 & 0 & 68921 & 122 & 82 & 79507 & 171 & 86 & 103823 & 138 & 0 \\ \hline
				48 & 50653 & 73 & 37 & 68921 & 121 & 0 & 79507 & 85 & 43 & 103823 & 46 & 0 \\ \hline
				49 & 50653 & 36 & 0 & 68921 & 122 & 82 & 79507 & 84 & 0 & 103823 & 46 & 0 \\ \hline
				50 & 50653 & 109 & 37 & 68921 & 81 & 41 & 79507 & 85 & 43 & 103823 & 187 & 141 \\ \hline
			\end{tabular}
	\end{table}}

\scriptsize{
	\begin{table}[p]
		\centering
		\caption{Proportions of Dimensions $n$ for Which Theorem~\ref{thm:P4-ideal} Proves Nonexistence, by Radius $r$ and Prime $p$.}
		\label{tab:appendix-ideal-prop}
		\begin{tabular}{c|c|c|c|c|c|c|c|c|c|c|c|c|c|c|c}
			\hline
			\diagbox[dir=SE]{$ r $}{proportion}{$ p $} & 3 & 5 & 7 & 11 & 13 & 17 & 19 & 23 & 29 & 31 & 37 & 41 & 43 & 47 & Total \\ \hline
			3 & 7.41\% & 0.00\% & 0.00\% & 0.00\% & 0.00\% & 0.00\% & 0.00\% & 0.00\% & 3.33\% & 0.00\% & 0.00\% & 0.00\% & 0.00\% & 0.00\% &10.49\%\\ \hline
			4 & 14.81\% & 0.00\% & 0.00\% & 0.00\% & 0.00\% & 0.00\% & 0.00\% & 0.00\% & 0.00\% & 0.00\% & 0.00\% & 0.00\% & 0.00\% & 0.00\%&14.81\%\\ \hline
			5 & 7.41\% & 3.20\% & 0.00\% & 0.00\% & 0.00\% & 0.00\% & 0.00\% & 0.00\% & 0.00\% & 0.00\% & 0.00\% & 0.00\% & 0.00\% & 0.00\%&10.37\%\\ \hline
			6 & 0.00\% & 7.20\% & 0.00\% & 0.00\% & 7.10\% & 0.00\% & 0.00\% & 0.00\% & 0.00\% & 0.00\% & 0.00\% & 0.00\% & 0.00\% & 0.00\%&13.79\%\\ \hline
			7 & 0.00\% & 21.60\% & 1.75\% & 0.00\% & 0.00\% & 0.00\% & 0.00\% & 0.00\% & 0.00\% & 0.00\% & 0.00\% & 0.00\% & 0.00\% & 0.00\%&22.97\%\\ \hline
			8 & 0.00\% & 3.20\% & 1.75\% & 0.00\% & 0.00\% & 5.54\% & 0.00\% & 0.00\% & 0.00\% & 0.00\% & 0.00\% & 0.00\% & 0.00\% & 0.00\%&10.16\%\\ \hline
			9 & 17.28\% & 3.20\% & 1.75\% & 0.00\% & 0.00\% & 0.00\% & 0.00\% & 0.00\% & 0.00\% & 0.00\% & 0.00\% & 0.00\% & 0.00\% & 0.00\%&21.33\%\\ \hline
			10 & 18.52\% & 6.40\% & 13.12\% & 0.00\% & 0.00\% & 0.00\% & 0.00\% & 0.00\% & 0.00\% & 0.00\% & 0.00\% & 0.00\% & 0.00\% & 0.00\%&33.74\%\\ \hline
			11 & 17.28\% & 6.40\% & 1.75\% & 0.75\% & 0.00\% & 0.00\% & 0.00\% & 0.00\% & 0.00\% & 0.00\% & 0.00\% & 0.00\% & 0.00\% & 0.00\%&24.50\%\\ \hline
			12 & 28.40\% & 27.20\% & 3.79\% & 1.50\% & 0.00\% & 0.00\% & 0.00\% & 0.00\% & 0.00\% & 0.00\% & 0.00\% & 0.00\% & 0.00\% & 0.00\%&50.60\%\\ \hline
			13 & 17.28\% & 10.40\% & 1.75\% & 1.50\% & 0.55\% & 0.00\% & 0.00\% & 0.00\% & 0.00\% & 0.00\% & 0.00\% & 0.00\% & 0.00\% & 0.00\%&28.67\%\\ \hline
			14 & 28.40\% & 6.40\% & 0.00\% & 1.50\% & 1.14\% & 0.00\% & 0.00\% & 0.00\% & 3.33\% & 0.00\% & 0.00\% & 0.00\% & 0.00\% & 0.00\%&36.91\%\\ \hline
			15 & 17.28\% & 4.00\% & 4.08\% & 1.50\% & 1.14\% & 0.00\% & 0.00\% & 0.00\% & 0.00\% & 0.00\% & 0.00\% & 0.00\% & 0.00\% & 0.00\%&25.83\%\\ \hline
			16 & 18.52\% & 8.00\% & 0.00\% & 8.64\% & 2.23\% & 0.00\% & 0.00\% & 0.00\% & 0.00\% & 0.00\% & 0.00\% & 0.00\% & 0.00\% & 0.00\%&33.04\%\\ \hline
			17 & 17.28\% & 12.00\% & 0.00\% & 1.50\% & 0.55\% & 0.33\% & 0.00\% & 0.00\% & 0.00\% & 0.00\% & 0.00\% & 0.00\% & 0.00\% & 0.00\%&28.93\%\\ \hline
			18 & 0.00\% & 4.00\% & 0.00\% & 1.50\% & 1.14\% & 0.67\% & 0.00\% & 0.00\% & 0.00\% & 0.00\% & 2.63\% & 0.00\% & 0.00\% & 0.00\%&9.59\%\\ \hline
			19 & 0.00\% & 4.00\% & 2.04\% & 1.50\% & 12.70\% & 1.34\% & 0.26\% & 0.00\% & 0.00\% & 0.00\% & 0.00\% & 0.00\% & 0.00\% & 0.00\%&20.43\%\\ \hline
			20 & 0.00\% & 0.00\% & 0.00\% & 2.33\% & 0.55\% & 0.67\% & 0.52\% & 0.00\% & 0.00\% & 0.00\% & 0.00\% & 4.76\% & 0.00\% & 0.00\%&8.59\%\\ \hline
			21 & 7.41\% & 0.00\% & 5.25\% & 0.75\% & 0.55\% & 0.33\% & 0.26\% & 0.00\% & 0.00\% & 0.00\% & 0.00\% & 0.00\% & 0.00\% & 0.00\%&13.91\%\\ \hline
			22 & 14.81\% & 16.00\% & 5.25\% & 0.00\% & 2.23\% & 1.02\% & 0.80\% & 0.00\% & 0.00\% & 0.00\% & 0.00\% & 0.00\% & 0.00\% & 0.00\%&34.91\%\\ \hline
			23 & 7.41\% & 4.00\% & 5.25\% & 0.00\% & 1.14\% & 0.65\% & 0.79\% & 0.54\% & 0.00\% & 0.00\% & 0.00\% & 0.00\% & 0.00\% & 0.00\%&18.37\%\\ \hline
			24 & 0.00\% & 0.00\% & 20.99\% & 0.00\% & 0.55\% & 1.00\% & 0.80\% & 0.18\% & 0.00\% & 0.00\% & 0.00\% & 0.00\% & 0.00\% & 0.00\%&22.97\%\\ \hline
			25 & 0.00\% & 13.44\% & 5.25\% & 0.00\% & 0.55\% & 10.20\% & 0.26\% & 0.18\% & 0.00\% & 0.00\% & 0.00\% & 0.00\% & 0.00\% & 0.00\%&27.07\%\\ \hline
			26 & 0.00\% & 11.68\% & 5.25\% & 0.00\% & 1.09\% & 1.00\% & 0.54\% & 0.75\% & 0.00\% & 0.00\% & 0.00\% & 0.00\% & 0.00\% & 0.00\%&19.11\%\\ \hline
			27 & 20.58\% & 22.40\% & 5.25\% & 0.00\% & 1.09\% & 1.34\% & 0.26\% & 0.36\% & 0.00\% & 0.00\% & 0.00\% & 2.38\% & 0.00\% & 0.00\%&44.72\%\\ \hline
			28 & 16.05\% & 10.88\% & 0.00\% & 0.00\% & 2.28\% & 1.36\% & 5.12\% & 0.73\% & 0.00\% & 0.00\% & 0.00\% & 0.00\% & 0.00\% & 0.00\%&32.07\%\\ \hline
			29 & 20.58\% & 13.44\% & 2.04\% & 0.00\% & 3.28\% & 0.33\% & 0.26\% & 0.18\% & 0.11\% & 0.00\% & 0.00\% & 0.00\% & 0.00\% & 0.00\%&35.43\%\\ \hline
			30 & 18.11\% & 16.64\% & 0.00\% & 0.00\% & 1.09\% & 0.67\% & 0.26\% & 1.11\% & 0.46\% & 0.00\% & 0.00\% & 0.00\% & 0.00\% & 0.00\%&34.16\%\\ \hline
			31 & 19.75\% & 13.92\% & 0.00\% & 1.65\% & 1.09\% & 1.34\% & 0.26\% & 0.18\% & 0.47\% & 0.10\% & 0.00\% & 0.00\% & 0.00\% & 0.00\%&34.38\%\\ \hline
			32 & 18.11\% & 24.48\% & 0.00\% & 0.00\% & 18.30\% & 0.33\% & 0.80\% & 0.73\% & 0.35\% & 0.10\% & 0.00\% & 0.00\% & 0.00\% & 0.00\%&50.63\%\\ \hline
			33 & 20.58\% & 14.72\% & 4.08\% & 0.75\% & 1.09\% & 0.33\% & 0.79\% & 0.54\% & 0.23\% & 0.10\% & 0.00\% & 0.00\% & 0.00\% & 0.00\%&37.48\%\\ \hline
			34 & 16.05\% & 16.64\% & 0.00\% & 1.50\% & 1.09\% & 0.65\% & 0.80\% & 4.25\% & 0.59\% & 0.51\% & 0.00\% & 2.38\% & 0.00\% & 0.00\%&37.88\%\\ \hline
			35 & 20.58\% & 19.84\% & 1.75\% & 1.50\% & 4.46\% & 1.34\% & 0.26\% & 0.54\% & 0.35\% & 0.40\% & 0.00\% & 0.00\% & 0.00\% & 0.00\%&42.83\%\\ \hline
			36 & 34.16\% & 17.76\% & 3.79\% & 1.50\% & 2.28\% & 1.65\% & 0.80\% & 0.73\% & 0.46\% & 0.20\% & 0.00\% & 0.00\% & 0.00\% & 0.00\%&51.76\%\\ \hline
			37 & 30.86\% & 26.56\% & 1.75\% & 1.50\% & 1.68\% & 1.34\% & 0.26\% & 0.18\% & 0.23\% & 0.20\% & 0.07\% & 0.00\% & 0.00\% & 0.00\%&52.79\%\\ \hline
			38 & 34.16\% & 16.96\% & 13.12\% & 8.64\% & 1.09\% & 0.65\% & 0.00\% & 1.11\% & 0.35\% & 0.71\% & 0.29\% & 0.00\% & 0.00\% & 0.00\%&58.39\%\\ \hline
			39 & 27.98\% & 19.84\% & 1.75\% & 1.50\% & 0.55\% & 2.04\% & 0.00\% & 0.18\% & 0.11\% & 0.10\% & 0.07\% & 0.00\% & 0.00\% & 0.00\%&45.82\%\\ \hline
			40 & 16.46\% & 16.64\% & 1.75\% & 1.50\% & 1.14\% & 1.65\% & 0.00\% & 0.73\% & 0.35\% & 0.20\% & 0.14\% & 0.00\% & 0.00\% & 0.00\%&35.41\%\\ \hline
			41 & 27.98\% & 13.92\% & 1.75\% & 1.50\% & 0.55\% & 1.30\% & 0.00\% & 0.36\% & 0.34\% & 0.30\% & 0.07\% & 0.06\% & 0.00\% & 0.00\%&41.78\%\\ \hline
			42 & 34.16\% & 29.28\% & 0.00\% & 2.33\% & 1.64\% & 14.86\% & 0.00\% & 0.75\% & 0.35\% & 0.41\% & 0.29\% & 0.18\% & 0.00\% & 0.00\%&62.66\%\\ \hline
			43 & 30.86\% & 14.72\% & 0.00\% & 0.75\% & 0.55\% & 1.30\% & 0.55\% & 0.18\% & 6.36\% & 0.31\% & 0.36\% & 0.18\% & 0.05\% & 0.00\%&47.09\%\\ \hline
			44 & 34.16\% & 16.64\% & 0.00\% & 0.00\% & 0.55\% & 1.30\% & 0.00\% & 0.37\% & 0.35\% & 0.10\% & 0.22\% & 0.29\% & 0.05\% & 0.00\%&46.86\%\\ \hline
			45 & 20.58\% & 13.44\% & 0.00\% & 0.83\% & 12.70\% & 1.69\% & 0.00\% & 0.54\% & 0.34\% & 0.10\% & 0.21\% & 0.24\% & 0.05\% & 0.00\%&42.35\%\\ \hline
			46 & 16.05\% & 10.88\% & 0.00\% & 0.00\% & 1.14\% & 0.65\% & 0.00\% & 0.00\% & 0.35\% & 3.17\% & 0.36\% & 0.18\% & 0.22\% & 0.00\%&29.63\%\\ \hline
			47 & 20.58\% & 22.40\% & 0.00\% & 0.00\% & 0.55\% & 0.65\% & 0.00\% & 0.38\% & 0.11\% & 0.10\% & 0.07\% & 0.18\% & 0.22\% & 0.13\%&39.82\%\\ \hline
			48 & 18.11\% & 10.88\% & 0.00\% & 0.00\% & 1.64\% & 1.99\% & 0.00\% & 0.00\% & 0.35\% & 0.20\% & 0.14\% & 0.18\% & 0.11\% & 0.04\%&30.36\%\\ \hline
			49 & 19.75\% & 13.44\% & 10.75\% & 0.00\% & 0.55\% & 1.00\% & 0.00\% & 0.38\% & 0.23\% & 0.31\% & 0.07\% & 0.18\% & 0.11\% & 0.04\%&39.75\%\\ \hline
			50 & 18.11\% & 26.88\% & 9.54\% & 0.00\% & 0.55\% & 0.65\% & 0.00\% & 0.00\% & 0.46\% & 0.41\% & 0.22\% & 0.12\% & 0.11\% & 0.18\%&47.27\%\\ \hline
		\end{tabular}
\end{table}}
	\restoregeometry
	
\end{document}